\RequirePackage{snapshot}
\documentclass[11pt,reqno]{amsart}
\usepackage{fixltx2e}                    
\usepackage[osf]{mathpazo}  
\usepackage[utf8]{inputenc}            
\usepackage{amsmath}                     
\usepackage{amssymb, latexsym, stmaryrd, amsthm, dsfont, amsfonts, amsbsy,amsthm, amsmath, mathrsfs}            
\usepackage{mathtools}                   
\usepackage{bm}                          
\usepackage{enumerate}                   
\usepackage{verbatim}                    
\usepackage{url}   
\usepackage{lscape}                      
\usepackage{microtype}  
\usepackage[all, knot]{xy}
        \xyoption{arc} 
        \xyoption{web}                 
\makeatletter                            
\def\MT@register@subst@font{\MT@exp@one@n\MT@in@clist\font@name\MT@font@list
 \ifMT@inlist@\else\xdef\MT@font@list{\MT@font@list\font@name,}\fi}
\makeatother

\usepackage[normal,sl,up,bf]{caption}	

\usepackage[pdftex,bookmarks,bookmarksnumbered,linktocpage,   
         colorlinks,linkcolor=blue,citecolor=blue]{hyperref}

\newcommand{\bit}{\begin{itemize}}    
\newcommand{\eit}{\end{itemize}}
\newcommand{\ben}{\begin{enumerate}}
\newcommand{\een}{\end{enumerate}}

\newcommand{\benroman}{\ben[\normalfont (i)]}  
\let\eroman\een

\newcommand{\bde}{\begin{description}}
\newcommand{\ede}{\end{description}}

\newcommand{\?}{\ensuremath{\mkern0.4\thinmuskip}}   
\let\models=\vDash                          

\let\leq=\leqslant
\let\nleq=\nleqslant
\let\geq=\geqslant

\let\epsilon=\varepsilon
\let\Lambda\varLambda
\let\Gamma\varGamma
\let\Delta\varDelta
\let\Lambda\varLambda
\let\Omega\varOmega
\let\Theta\varTheta
\let\Xi\varXi
\let\Pi\varPi
\let\Sigma\varSigma

\let\bs=\boldsymbol
\let\class=\mathsf                              
\let\oper=\mathbb                               

\bmdefine{\A}{A}                                
\bmdefine{\2}{2}
\bmdefine{\B}{B}
\bmdefine{\D}{D}
\bmdefine{\M}{M}                                
\bmdefine{\LLL}{L}                              
\bmdefine{\Fm}{Fm}                              
\bmdefine{\zerou}{[0{,}1]}  
\bmdefine{\T}{T}                                

\newcommand{\VVV}{\oper{V}}                     

\newcommand{\HHH}{\oper{H}}

\newcommand{\PPU}{\oper{P}_{\!\textsc{u}}^{}}

\newcommand{\SSS}{\oper{S}}
\newcommand{\III}{\oper{I}}

\bmdefine{\boldstar}{\mathchoice{\textstyle*}{\textstyle*}{\textstyle*}{\scriptstyle*}}

\bmdefine{\btau}{\tau}                                  
\bmdefine{\brho}{\rho}                                  

\bmdefine{\leibniz}{\Omega}        

\bmdefine{\frege}{\Lambda}         

\makeatletter
\newcommand{\tarskidsp}{\mathord%
   {\m@th\raisebox{0pt}[0pt][0pt]{$\stackrel%
   {\raisebox{-2.7pt}[0ex][0pt]{$\displaystyle \,\?\thicksim$}}%
   {\displaystyle\leibniz}$}}}
\newcommand{\tarskitxt}{\mathord%
   {\m@th\raisebox{0pt}[0pt][0pt]{$\stackrel%
   {\raisebox{-2.7pt}[0ex][0pt]{$\,\?\thicksim$}}{\displaystyle\leibniz}$}}}
\newcommand{\tarskiscr}{\mathord%
   {{\m@th\raisebox{0pt}[0pt][0pt]{$\stackrel%
   {\raisebox{-2.4pt}[0ex][0pt]{$\scriptstyle \,\?\thicksim$}}%
   {\scriptstyle\leibniz}$}}}}
\newcommand{\tarskiscrscr}{\mathord%
   {{\m@th\raisebox{0pt}[0pt][0pt]{$\stackrel%
   {\raisebox{-2pt}[0ex][0pt]{$\scriptscriptstyle \,\?\thicksim$}}%
   {\scriptscriptstyle\leibniz}$}}}}
\newcommand{\tarski}{\@ifnextchar ^ %
   {\mathchoice{\tarskidsp\kern-.07em}{\tarskitxt\kern-.07em}%
   {\tarskiscr\kern-.07em}{\tarskiscrscr\kern-.07em}}%
   {\mathchoice{\tarskidsp}{\tarskitxt}{\tarskiscr}{\tarskiscrscr}}}
\makeatother

\theoremstyle{theorem}
\newtheorem{Theorem}{Theorem}[section]
\newtheorem{Lemma}[Theorem]{Lemma}
\newtheorem{Corollary}[Theorem]{Corollary}
\newtheorem{Proposition}[Theorem]{Proposition}

\theoremstyle{definition}

\theoremstyle{remark}

\newtheorem{Remark}[Theorem]{Remark}

\newcommand{\C}{\boldsymbol{C}} 

\newsavebox{\XnLines}
\savebox{\XnLines}(72,24)[bl]{

\put(0,0){\line(3,5){14}}
\put(0,24){\line(3,-5){14}}
\put(14,0){\line(0,1){24}}

\put(0,0){\line(6,5){28}}
\put(0,24){\line(6,-5){28}}
\put(28,0){\line(0,1){24}}

\put(0,0){\line(3,1){72}}
\put(0,24){\line(3,-1){72}}
\put(72,0){\line(0,1){24}}

}

\newsavebox{\XnPic}
\savebox{\XnPic}(72,24)[bl]{%

\put(0,0){\usebox{\XnLines}}

\put(0,0){\circle*{4}}
\put(0,24){\circle*{4}}
\put(14,24){\circle*{4}}
\put(14,0){\circle*{4}}

\put(28,24){\circle*{4}}
\put(28,0){\circle*{4}}

\put(72,24){\circle*{4}}
\put(72,0){\circle*{4}}

\put(44,0){\circle*{2}}
\put(50,0){\circle*{2}}
\put(56,0){\circle*{2}}
\put(44,24){\circle*{2}}
\put(50,24){\circle*{2}}
\put(56,24){\circle*{2}}
}

\begin{document}
\title[Epimorphism surjectivity in varieties of Heyting algebras]{Epimorphism surjectivity in varieties of Heyting algebras}

\author{T. Moraschini}
\address{Department of Philosophy, University of Barcelona, Carrer de Montalegre $6$, $08001$, Barcelona, Spain\newline {\rm and } \newline Institute of Computer Science, Academy of Sciences of Czech Republic, Pod Vod\'arenskou v\v{e}\v{z}\'{i} $271/2$, $182$ $07$ Prague $8$, Czech Republic}\email{tommaso.moraschini@ub.edu}

\author{J. J. Wannenburg}
\address{Department of Mathematics and Applied Mathematics, University of Pretoria, Private Bag X$20$, Hatfield, Pretoria $0028$, and DST-NRF Centre of Excellence in Mathematical and Statistical Sciences (CoE-MaSS), South Africa}
\email{jamie.wannenburg@up.ac.za}
\renewcommand{\thefootnote}{\fnsymbol{footnote}} 
\footnotetext{\emph{Key words.} Epimorphism, Heyting algebra, Esakia space, intuitionistic logic, intermediate logic, Beth definability.}     
\footnotetext{\emph{2010 Mathematics Subject Classification.} $03$B$55$, $06$D$20$, $18$A$20$, $06$D$50$, $03$G$27$, $03$G$10$.}     
\renewcommand{\thefootnote}{\arabic{footnote}} 
\date{\today}

\maketitle


\begin{abstract}
It was shown recently that epimorphisms need not be surjective in a variety $\mathsf{K}$ of Heyting algebras, but only
one counter-example was exhibited in the literature until now.
Here, a continuum of such examples is identified, viz.\ the variety generated by the Rieger-Nishimura lattice, and all of its (locally finite) subvarieties that contain the original counter-example $\mathsf{K}$.  It is known that, whenever a variety of Heyting algebras has finite depth, then it has surjective epimorphisms.  In contrast, we show that for every integer $n\geqslant 2$, the variety of all Heyting algebras of width at most $n$ has a non-surjective epimorphism.
Within the so-called Kuznetsov-Ger\v{c}iu variety (i.e., the variety generated by finite linear sums of one-generated Heyting algebras), we describe exactly the subvarieties that have surjective epimorphisms.
This yields new positive examples, 
%
%
and an alternative proof of epimorphism surjectivity for all varieties of G\"{o}del algebras.  The results
settle natural questions about Beth-style definability for a range of intermediate logics.
\end{abstract}

\section{Introduction}


A morphism $f \colon \A \to \B$ in a category is called an \emph{epimorphism}
\cite{Bacsich74,Cam18,Isb66,KMPT83}
provided that it is right-cancellative, i.e., for every pair of morphisms $g, h \colon \B \to \C$ (in the same category),
\[
\text{if }g \circ f = h \circ f \text{, then }g= h.
\]
We regard any variety of algebras as a concrete category whose
morphisms are the algebraic homomorphisms between its members.
In such categories all surjective morphisms are
epimorphisms, but
the converse is
not
true in general. Accordingly, when all epimorphisms are surjective in a variety $\class{K}$ of algebras, we say that $\class{K}$ has the \emph{epimorphism surjectivity (ES) property}.  

This requirement can be phrased both in categorical and model theoretic terms. On the one hand, an epimorphism $f \colon \A \to \B$ is said to be \emph{regular} if it is the co-equalizer of a pair of morphisms $g, h \colon \C \to \A$. 
In a variety it can be shown that the co-equalizer of a pair of homomorphisms $g,h \colon \C \to \A$ is the canonical surjection $f \colon \A \to \A/\theta$, where $\theta$ is the smallest congruence of $\A$ containing the set $\{\langle g(c),h(c) \rangle \in A \times A \colon c \in C\}$.
It follows that every regular epimorphism is a surjective homomorphism in a variety of algebras.
The converse is also true, as every surjective homomorphism $f \colon \A \to \B$ can be seen as the co-equalizer of a pair of maps $p_{1}, p_{2} \colon \textup{Ker}(f) \to \A$, where $\textup{Ker}(f)$ is the kernel of $f$ seen as a subalgebra of $\A \times \A$ and $p_{1}$ and $p_{2}$ are the projection maps, respectively, on the first and on the second component. Consequently, in varieties of algebras, regular epimorphisms coincide with surjective homomorphisms, and a variety has the ES property if and only if epimorphisms are regular in it. 

On the other hand, from the point of view of model theory, the ES property is related to definability. More precisely, a variety $\class{K}$ lacks the ES property when there are two algebras $\A, \B \in \class{K}$ such that $\A$ is a proper subalgebra of $\B$, and for every $b \in B \smallsetminus A$ there exist a primitive positive formula $\varphi(\vec{x}, y)$ and a tuple $\vec{a} \in A$ such that $\B \vDash \varphi(\vec{a}, b)$ and $\varphi(\vec{x}, y)$ defines a partial function in $\class{K}$ \cite[Thm.\ 3.1]{Cam18}.


The ES property need not be inherited by subvarieties. For example, the variety $\class{L}$ of all lattices has the ES property, whence an embedding $f$ of the three-element chain into
the four-element diamond is not an
epimorphism in $\class{L}$.  Indeed,
the
diamond can be embedded into the five-element non-distributive modular lattice 
in two distinct ways that agree when composed with $f$.
Nevertheless, $f$ is a (non-surjective) epimorphism in the smaller variety of distributive lattices.
As this example suggests,
it is
often difficult to establish whether epimorphisms are surjective in a given variety.

The failure of the ES property for distributive lattices can be explained in logical terms by the observation that complements are implicitly but not explicitly definable---i.e., when complements exist in a distributive lattice, they are uniquely determined, but
no unary term defines them explicitly.  This instantiates a general result:
the algebraic counterpart $\mathsf{K}$ of an algebraizable logic $\vdash$ \cite{BP89} has the ES property if and only if $\vdash$ has
the \emph{infinite (deductive) Beth (definability) property}---i.e.,
all implicit definitions in $\vdash$ can be made explicit \cite{BH06,Hoo01}.\footnote{The focus on algebraizable logics is not restrictive, because
every (pre)variety is categorically equivalent to one that algebraizes some sentential logic \cite{MorRaf18-prevarieties}, and the ES property persists under any category equivalence between (pre)varieties.} The
terminology comes from Beth's definability theorem for classical predicate logic \cite{Bet53}.

More precisely, recall that every algebraizable logic $\vdash$ has a set of formulas $\Delta(x, y)$ which behaves globally as an equivalence connective \cite{BP89,Cze01,Fon16,FJa09, FJP03}.
Consider two
disjoint sets $X$ and $Z$ of variables, with $X \ne \emptyset$, and a set $\Gamma$ of formulas over $X \cup Z$.  We say that $Z$ is \emph{defined implicitly} in terms of $X$ by means of $\Gamma$ in $\vdash$ if for all $z \in Z$ 
\[
\Gamma \cup \sigma[\Gamma] \vdash \Delta(z, \sigma(z))
\]
for every substitution $\sigma$ such that $\sigma(x) = x$ for all $x \in X$. On the other hand, $Z$ is said to be \emph{defined explicitly} in terms of $X$ by means of $\Gamma$ in $\vdash$ when, for every $z \in Z$, there exists a formula $\varphi_{z}$ over $X$ such that
\[
\Gamma \vdash \Delta(z, \varphi_{z}).
\]
Then the infinite Beth property postulates the equivalence of implicit and explicit definability in $\vdash$
(for all $X,Z,\Gamma$ as above).\footnote{Observe that $X, Z$ and $\Gamma$ may be arbitrarily large sets, and that $\sigma$ can map $X \cup Z$ to sets of variables other than $X \cup Z$. For logics satisfying suitable (possibly infinitary) versions of compactness, some cardinal bounds can be imposed on $X,Z$, and on the codomain of $\sigma$, as is shown in \cite{MRW3} (see also \cite{Bacsich74,Isb66}).
} The \emph{finite Beth property} does the same for the case where $Z$ is a \emph{finite} set.

It turns out that an algebraizable logic has the finite Beth property if and only if its algebraic counterpart has the
\emph{weak ES property} \cite{BH06,Hoo01},\footnote{This result has an antecedent, due to N\'emeti, in \cite[Thm.\ 5.6.10]{HMT7185}.}
which means that the ``almost onto'' epimorphisms between its members are surjective. Here, a homomorphism $f \colon \A \to \B$ is said to be \emph{almost onto} if there is a \emph{finite} subset $C \subseteq B$ such that $f[A] \cup C$ generates $\B$.

The foregoing facts motivate the study of epimorphisms in \emph{varieties of Heyting algebras}, i.e., in the algebraic counterparts of intermediate logics.
A well-known result of Kreisel states that all intermediate logics have the finite Beth property, whence all \emph{all} varieties of Heyting algebras have the \emph{weak} ES property \cite{Kre60}. In contrast, a classical 
result of Maksimova states that only \emph{finitely many} such varieties possess a \emph{strong} variant of the ES property \cite{GM05,Mak99ail,Mak00,Mak03}, i.e., satisfy the requirement that if $f \colon \A \to \B$ is a homomorphism and $b \in B \smallsetminus A$, then there is a pair of homomorphisms $g, h \colon \B \to \C$ such that $g(b) \ne h(b)$ and $g \circ f = h \circ f$. Varieties of Heyting algebras with the strong ES property include the respective varieties of all Heyting algebras, all G\"odel algebras \cite{Dm59}, and all Boolean algebras. 
Notably, the strong ES property admits a categorical formulation as the demand that all monomorphisms are regular, i.e., are the equalizer of some parallel pair of morphisms.\footnote{Recall that, in general, a monomorphism is a left-cancellative morphism. However, in varieties, monomorphisms coincide with algebraic embeddings. Therefore, as equalizers are always monomorphisms, when applied to varieties, the demand that all monomorphisms are regular amounts to the requirement that embeddings and equalizers coincide.}

Very little is known, however, about the (unqualified) ES property for varieties of Heyting algebras, despite the fact that it is algebraically natural. On the positive side, it was shown recently that the ES property holds for all varieties of Heyting algebras
that have
finite
depth \cite[Thm.\ 5.3]{BMR17}. This yields a continuum of examples \emph{with} the ES property.
On the other hand, in the same paper \cite[Cor.\ 6.2]{BMR17} it was shown that even a locally finite variety
of Heyting algebras need not have the ES property. (The counter-example
confirmed Blok and Hoogland's conjecture \cite{BH06} that the weak ES property really is strictly weaker than the ES property.) This raises the question: how rare are varieties of Heyting algebras \emph{without} the ES property?

In this work we establish that the ES property fails for the variety of all Heyting algebras of
width at most $n+2$, where $n$ is any natural number (Corollary \ref{Cor:width-n}).
We also disprove the ES property for the variety generated by the Rieger-Nishimura lattice, and for a continuum of its locally finite subvarieties (Corollary \ref{Cor:RN-not-ES} and Theorem \ref{Thm:continuum}).
Finally,
within the \emph{Kuznetsov-Ger\v{c}iu variety} $\class{KG}$,
i.e., the variety generated by finite linear sums of one-generated Heyting algebras,
we classify the subvarieties with the ES property
(Theorem \ref{Thm:ES-in-KG}).  As $\class{KG}$ contains all G\"{o}del algebras,
this gives a new explanation of the surjectivity of epimorphisms in all varieties of G\"{o}del algebras
(which was first shown in \cite[Cor.\ 5.7]{BMR17}).

\section{Esakia duality}

A \emph{Heyting algebra} is an algebra $\A = \langle A; \land, \lor, \to, 0, 1 \rangle$ which comprises a bounded lattice $\langle A; \land, \lor, 0, 1 \rangle$, and a binary operation $\to$ such that for every $a, b, c \in A$,
\[
a \land b \leq c \Longleftrightarrow a \leq b \to c.
\]
It follows that Heyting algebras are distributive lattices. Remarkably, a Heyting algebra is uniquely determined by its lattice reduct. The class of all Heyting algebras forms a variety, which we denote by $\class{HA}$. 

A fundamental tool for investigating Heyting algebras is Esakia duality \cite{Esa74,Esa85,Celani2014}, which we proceed to review. Given a poset $\langle X; \leq \rangle$, we call any upward closed subset of $X$ an \emph{upset}. Similarly, any downward closed subset is called a \emph{downset}. For any set $U \subseteq X$, the smallest upset and downset containing $U$ are denoted respectively by ${\uparrow} U$ and ${\downarrow} U$. In case $U = \{ x \}$, we shall write ${\uparrow} x$ and ${\downarrow} x$ instead of ${\uparrow} \{ x \}$ and ${\downarrow} \{ x\}$, respectively. Then an \emph{Esakia space} $\boldsymbol{X} = \langle X; \tau, \leq \rangle$ comprises a Stone space $\langle X; \tau \rangle$ (i.e., a compact Hausdorff space in which every open set is a union of clopen sets) and a poset $\langle X; \leq \rangle$ such that
\benroman
\item ${\uparrow} x$ is closed for all $x \in X$, and 
\item ${\downarrow} U$ is clopen, for every clopen $U \subseteq X$. 
\eroman
Observe that the topology of finite Esakia spaces is necessarily discrete (because they are Hausdorff), and that finite posets endowed with the discrete topology are Esakia spaces. Moreover, every Esakia space $\bs{X}$ satisfies the \emph{Priestley separation axiom} \cite{Pr70,Pr72}, stating that for every $x, y  \in X$ such that $x \nleq y$, there exists a clopen upset $U$ such that $x \in U$ and $y \notin U$.

For Esakia spaces $\boldsymbol{X}$ and $\boldsymbol{Y}$, an \emph{Esakia morphism} $f \colon \boldsymbol{X} \to \boldsymbol{Y}$ is a continuous order-preserving map $f \colon X \to Y$ such that for all $x \in X$,
\begin{equation}\label{Eq:Esakia-morphism}
\text{if } f(x) \leq y \in Y, \text{ then } y = f(z) \text{ for some } z \geq x.
\end{equation}
Esakia spaces form a category, which we denote by $\class{ESP}$, in which the morphisms are Esakia morphisms. 
Note that isomorphisms of $\class{ESP}$ are exactly bijective Esakia morphisms, because Esakia spaces are compact Hausdorff.

The relation between Heyting algebras and Esakia spaces is as follows: 
\begin{Theorem}[L.\ Esakia {\cite[Thm.\ 3, p.\ 149]{Esa74}}]\label{Thm:Esakia}
The categories $\class{HA}$ and $\class{ESP}$ are dually equivalent.
\end{Theorem}

The dual equivalence functors are defined as follows. Given a Heyting algebra $\A$, we denote the set of its (non-empty proper) prime filters by $\textup{Pr}\A$. For every $a \in A$, set 
\begin{equation}\label{Eq:gamma}
\gamma^{\A}(a) \coloneqq \{ F \in \textup{Pr}\A \colon a \in F \},
\end{equation}
and consider also its complement $\gamma^{\A}(a)^{c} \coloneqq \{F \in \textup{Pr}\A \colon a \notin F \}$.
It turns out that the structure $\A_{\ast} \coloneqq \langle \textup{Pr}\A; \tau, \subseteq \rangle$ is an Esakia space, where $\tau$ is the topology on $\textup{Pr}\A$ with subbasis $\{ \gamma^{\A}(a) \colon a \in A \} \cup \{ \gamma^{\A}(a)^{c} \colon a \in A \}$. Moreover, for every Heyting algebra homomorphism $f \colon \A \to \B$, let $f_{\ast} \colon \B_{\ast} \to \A_{\ast}$ be the Esakia morphism defined by the rule $F \mapsto f^{-1}[F]$.

Conversely, let $\boldsymbol{X}$ be an Esakia space. We denote by $\textup{Cu}\boldsymbol{X}$ the set of clopen upsets of $\boldsymbol{X}$. Then the structure $\boldsymbol{X}^{\ast} \coloneqq \langle \textup{Cu}\boldsymbol{X}; \cap, \cup, \to, \emptyset, X \rangle$, where $U \to V \coloneqq X \smallsetminus {\downarrow} (U \smallsetminus V)$, is a Heyting algebra. Moreover, for every Esakia morphism $f \colon \boldsymbol{X} \to \boldsymbol{Y}$, let $f^{\ast} \colon \boldsymbol{Y}^{\ast} \to \boldsymbol{X}^{\ast}$ be the homomorphism of Heyting algebras given by the rule $U \mapsto f^{-1}[U]$. 

The dual equivalence in Theorem \ref{Thm:Esakia} is witnessed by the pair of contravariant functors
\begin{equation}\label{Eq:functors}
(-)_{\ast} \colon \class{HA} \longleftrightarrow \class{ESP} \colon (-)^{\ast}.
\end{equation}
Given a variety $\mathsf{K}$ of Heyting algebras, we denote by $\class{K}_{\ast}$ the full subcategory of $\class{ESP}$ 
whose class of objects is the isomorphic closure of $\{ \A_{\ast} \colon \A \in \mathsf{K} \}$. It is clear that the functors in (\ref{Eq:functors}) restrict to a dual equivalence between $\mathsf{K}$ and $\mathsf{K}_{\ast}$.

Let $\boldsymbol{X}$ be an Esakia space. An \emph{Esakia subspace} (E-subspace for short) of $\boldsymbol{X}$ is a closed upset of $\boldsymbol{X}$, equipped with the subspace topology and the restriction of the order. A \emph{correct partition} on $\boldsymbol{X}$ (sometimes called an \emph{Esakia relation} or \emph{bisimulation equivalence} in the literature) is an equivalence relation $R$ on $X$ such that for every $x, y, z \in X$,
\benroman
\item if $\langle x,y \rangle \in R$ and $x \leq z$, then $\langle z,w \rangle \in R$ for some $w \geq y$, and
\item if $\langle x, y \rangle \notin R$, then there is a clopen $U$, 
such that $x \in U$ and $y \notin U$, which moreover is a union of equivalence classes of $R$.
\eroman
In this case, we denote by $\bs{X} / R$ the Esakia space consisting of the quotient space of $\bs{X}$ with respect to $R$, equipped with the partial order $\leq^{\bs{X}/R}$ defined as follows: for every $x, y \in X$,
\begin{align*}
x/ R \leq^{\bs{X}/R} y/R \Longleftrightarrow& \text{ there are }x', y' \in X \text{ such that}\\
&\text{ }\langle x, x' \rangle, \langle y, y' \rangle \in R\text{ and }x' \leq^{\bs{X}}y'.
\end{align*}
The map $x \mapsto x/R$ is an Esakia morphism from $\bs{X}$ to $\bs{X}/R$, and 
for every Esakia morphism $f \colon \bs{X} \rightarrow \bs{Y}$, the kernel of $f$ is a correct partition on $\bs{X}$. If, moreover, $f$ is surjective, then there is an $\class{ESP}$-isomorphism $i \colon \bs{X}/ \ker f \cong \bs{Y}$, such that $i \circ q = f$.

An algebra $\A$ is \emph{finitely subdirectly irreducible}, FSI for short, when the identity relation is meet-irreducible in the congruence lattice of $\A$ \cite{Ber12,BS81}.

\begin{Lemma}\label{Lem:correspondences}
 Let $\A$ be a Heyting algebra, and $\boldsymbol{X}$ an Esakia space.
\benroman
\item\label{Lem:correspondences:FSI} $\A$ is FSI if and only if its top element is \emph{prime} (i.e., if $x \vee y = 1$ then $x = 1$ or $y = 1$) or, equivalently, the poset underlying $\A_{\ast}$ is \emph{rooted} (i.e., has a least element).
\item \label{Lem:correspondences:InjectiveSurjective} A homomorphism $h$ between Heyting algebras is injective iff $h_{\ast}$ is surjective. Also, $h$ is surjective iff $h_{\ast}$ is injective.
\item\label{Lem:correspondences:FactorAndSubspaces} There is a dual lattice isomorphism $\sigma$ from the congruence lattice of $\A$ to the lattice of E-subspaces of $\A_{\ast}$, such that $(\A/\theta)_* \cong \sigma(\theta)$ for any congruence $\theta$ of $\A$, and for any E-subspace $\bs{Y}$ of $\A_{\ast}$, we have $\bs{Y}^{\ast} \cong \A/\sigma^{-1}(\bs{Y})$.
\item\label{Lem:correspondences:SubalgebrasAndPartitions} There is a dual lattice isomorphism $\rho$ from the lattice of subalgebras of $\A$ to that of correct partitions on $\A_{\ast}$, such that if $\B$ is a subalgebra of $\A$ then $\B_{\ast} \cong \A_{\ast}/\rho(\B)$, and if $R$ is a correct partition on $\A_{\ast}$ then $(\A_{\ast}/R)^{\ast} \cong \rho^{-1}(R)$.
\item\label{Lem:correspondences:SubspacesOfPartitions} Let $R$ be a correct partition on $\boldsymbol{X}$. If $\boldsymbol{Y}$ is an E-subspace of $\boldsymbol{X}$, then $R \cap Y^{2}$ is a correct partition on $\boldsymbol{Y}$. 
\item\label{Lem:correspondences:ImagesAreSubspaces} Images of Esakia morphisms are E-subspaces, and restrictions of Esakia morphisms to E-subspaces are still Esakia morphisms.
\item\label{Lem:correspondences:Chains} Every chain in $\boldsymbol{X}$ has an infimum and a supremum. Moreover, infima and suprema of chains are preserved by Esakia morphisms.
\eroman
\end{Lemma}
The statement of (\ref{Lem:correspondences:FSI}) is well-known (see for instance \cite[Thm.\ 2.9]{BeBe08}). The correspondences of (\ref{Lem:correspondences:InjectiveSurjective}), (\ref{Lem:correspondences:FactorAndSubspaces}) and (\ref{Lem:correspondences:SubalgebrasAndPartitions}) were established in \cite{Esa74} (alternatively, examine \cite[Lem.~3.4]{BMR17}).
The isomorphism $\sigma$ in (\ref{Lem:correspondences:FactorAndSubspaces}) is the map that sends a congruence $\theta$ of $\A$ to the E-subspace of $\A_{\ast}$ with universe
\[
\{ F \in \textup{Pr}\A \colon F \text{ is a union of equivalence classes of }\theta \}.
\]
Moreover, the isomorphism $\rho$ in (\ref{Lem:correspondences:SubalgebrasAndPartitions}) is the map that sends a subalgebra $\A$ of $\B$ (in symbols $\A \leq \B$) to the correct partition $R_{\A}$ on $\B_{\ast}$ such that
\begin{equation}\label{Eq:correct-partition}
\langle F, G \rangle \in R_{\A} \Longleftrightarrow F \cap A = G \cap A
\end{equation}
for every $F, G \in \textup{Pr}\B$. 

To prove (\ref{Lem:correspondences:SubspacesOfPartitions}), notice that $R \cap Y^2$ is the kernel of the Esakia morphism obtained by composing the canonical surjection from $\bs{X}$ to $\bs{X}/R$ with the inclusion map from $\bs{Y}$ to $\bs{X}$. 
Moreover, the first part of (\ref{Lem:correspondences:ImagesAreSubspaces}) is a consequence of (\ref{Eq:Esakia-morphism}) and of the fact that continuous functions from compact spaces to Hausdorff spaces are closed. The second part of (\ref{Lem:correspondences:ImagesAreSubspaces}) is immediate. Finally, (\ref{Lem:correspondences:Chains}) is a consequence of Theorem \ref{Thm:Esakia}, together with the observation that unions and intersections of chains of prime filters are still prime filters (and that these unions and intersections are preserved by inverse images of homomorphisms).

Note that, owing to (\ref{Lem:correspondences:FactorAndSubspaces}) and (\ref{Lem:correspondences:SubalgebrasAndPartitions}), if $\class{K}$ is a variety of Heyting algebras then $\class{K}_{\ast}$ is closed under taking E-subspaces and quotients by correct partitions.

\section{Epimorphism surjectivity}

Let $\class{K}$ be a variety of algebras and $\B \in \mathsf{K}$. A subalgebra $\A \leq \B$ is $\class{K}$-\emph{epic} if for every pair of morphisms $g, h \colon \B \to \C$ in $\mathsf{K}$,
\[
\text{if }g \! \upharpoonleft_{A} = h \! \upharpoonleft_{A}\text{, then }g=h.
\]

\begin{Lemma}\label{Lem:epic-subalgebras}
A variety $\mathsf{K}$ has the ES property if and only if 
no algebra in $\class{K}$ has a proper $\class{K}$-epic subalgebra.
\end{Lemma}

\begin{proof}
Observe that if there is a non-surjective epimorphism $f \colon \A \to \B$ in $\mathsf{K}$, then $f[\A]$ is a proper $\mathsf{K}$-epic subalgebra of $\B$. Conversely, if $\A \leq \B$ is a proper $\mathsf{K}$-epic subalgebra of $\B$, then the inclusion map $\A \hookrightarrow \B$ is a non-surjective epimorphism in $\mathsf{K}$.
\end{proof}

The next result is a topological reformulation of the above result in the special case of Heyting algebras, which is essentially a consequence of the correspondence between subalgebras and correct partitions (see (\ref{Lem:correspondences:InjectiveSurjective}) and (\ref{Lem:correspondences:SubalgebrasAndPartitions}) of Lemma \ref{Lem:correspondences}). Also, because of the dual equivalence, the dual of every epimorphism is a \emph{monomorphism} (i.e., a morphism $f$ such that, for any morphisms $g$ and $h$, if $f \circ g = f \circ h$ then $g = h$).



\begin{Lemma}\label{Lem:epic-subalgebras-dual}
A variety $\mathsf{K}$ of Heyting algebras lacks the ES property if and only if there is an Esakia space $\boldsymbol{X} \in \class{K}_{\ast}$ with a correct partition $R$ different from the identity relation such that for every $\boldsymbol{Y} \in \class{K}_{\ast}$ and every pair of Esakia morphisms $g, h \colon \boldsymbol{Y} \to \boldsymbol{X}$, if $\langle g(y), h(y) \rangle \in R$ for every $y \in Y$, then $g= h$. 
\end{Lemma}


A variety is said to be \emph{arithmetical} when it is both congruence permutable and congruence distributive \cite{Ber12,BS81}. Given a class of algebras $\class{K}$, we denote by $\class{K}_{\textup{FSI}}$ the class of all FSI members of $\class{K}$. We say $\class{K}$ is a \emph{universal class} when $\class{K}$ is axiomatized by universal first-order sentences. 

\begin{Theorem}[M.\ Campercholi {\cite[Thm.\ 6.8]{Cam18}}]\label{Thm:Campercholi}
Let $\mathsf{K}$ be an arithmetical variety such that $\class{K}_{\textup{FSI}}$ is a universal class.\ Then $\mathsf{K}$ has the ES property if and only if its FSI members lack proper $\mathsf{K}$-epic subalgebras.
\end{Theorem}

Notice from Lemma \ref{Lem:correspondences}(\ref{Lem:correspondences:FSI}) that FSI Heyting algebras form a universal class. Moreover, it is well-known that the variety of Heyting algebras is arithmetical \cite{GJKO07}. As a consequence we can instantiate Theorem \ref{Thm:Campercholi} as follows: 

\begin{Corollary}\label{Cor:FSI}
A variety $\mathsf{K}$ of Heyting algebras has surjective epimorphisms if and only if its FSI members lack proper $\mathsf{K}$-epic subalgebras.
\end{Corollary}

To illustrate the power of Corollary \ref{Cor:FSI}, we shall exhibit a new elementary proof of the fact that finitely generated varieties of Heyting algebras (i.e., varieties generated by a finite algebra) have surjective epimorphisms \cite[Cor.\ 5.5]{BMR17}. 
To this end, 
we denote by 
$\HHH, \SSS, \PPU$ and $\VVV$ 
the respective class operators for homomorphic images, subalgebras, ultraproducts and varietal generation. Recall \emph{J\'onsson's lemma}, which states that if $\mathsf{K}$ is a class of similar algebras and $\VVV(\mathsf{K})$ is congruence distributive, then $\VVV(\mathsf{K})_{\textup{FSI}} \subseteq \HHH\SSS\PPU(\mathsf{K})$, see \cite[Lem.\ 5.9]{Ber12} or \cite{Jon67,Jon95}. 
In particular, if $\class{K}$ is a finite set of finite algebras, then $\VVV(\class{K})_{\textup{FSI}} \subseteq \HHH\SSS(\class{K})$.

\begin{Proposition}\label{Prp:FG-weak-ES-implies-ES}
If a finitely generated variety 
of Heyting algebras 
has the weak ES property then it has the (unqualified) ES property.
\end{Proposition}
\begin{proof}
Let $\class{K}$ be a finitely generated variety of Heyting algebras without the ES property.
By Corollary \ref{Cor:FSI}, there is an FSI algebra $\B \in \mathsf{K}$ with a proper $\class{K}$-epic subalgebra $\A$. Thus the inclusion map $\A \hookrightarrow \B$ is a non-surjective epimorphism. By J\'onsson's lemma and the fact that $\mathsf{K}$ is finitely generated, we know that $\B$ is finite. Therefore the map $\A \hookrightarrow \B$ is almost onto. Thus $\mathsf{K}$ has a non-surjective almost onto epimorphism, i.e., it lacks the \emph{weak} ES property.
\end{proof}
More generally,
by appealing to \cite[Cor.~6.5]{Cam18} instead of Corollary \ref{Cor:FSI}, 
the proof of Proposition~\ref{Prp:FG-weak-ES-implies-ES} can be generalized to show that
the weak ES property entails the ES property for
every finitely generated variety of algebras with a majority term (e.g.\ one generated by a finite lattice-based algebra).
\begin{Theorem}[G.\ Kreisel {\cite{Kre60}; also see \cite[Sec.\ 12]{GR15}}]\label{Thm:Kreisel}
All varieties of Heyting algebras have the weak ES property.
\end{Theorem}
Proposition~\ref{Prp:FG-weak-ES-implies-ES} together with Theorem~\ref{Thm:Kreisel} yields the following: 

\begin{Proposition}[{\cite[Cor.\ 5.5]{BMR17}}]\label{Prp:finitely-generated}
Every finitely generated variety 
of Heyting algebras 
has surjective epimorphisms.
\end{Proposition}

\section{Depth and width in Heyting algebras}

Let $0 < n \in \omega$. A poset $\mathbb{X} = \langle X; \leq \rangle$ is said to have \emph{depth at most} $n$ if it does not contain any chain of $n+1$ elements. Similarly, $\mathbb{X}$ is said to have \emph{width at most} $n$ if ${\uparrow} x$ does not contain any antichain of $n+1$ elements, for every $x \in X$. A Heyting algebra $\A$ has \emph{depth} (resp.\ \emph{width}) \emph{at most} $n$, when the poset underlying its dual space $\A_{\ast}$ has depth (resp.\ width) at most $n$. 


For $0 < n \in \omega$, let $\class{D}_{n}$ and $\class{W}_{n}$ be, respectively, the classes of Heyting algebras with depth and width at most $n$. It follows that $\class{D}_{1}$ is the variety of Boolean algebras, while $\class{W}_{1}$ is the variety of G\"odel algebras (i.e., of subdirect products of totally ordered Heyting algebras). Both of them are known to have a strong variant of the ES property \cite{GM05,Mak03,Mak99ail,Mak00}. The following result is well known (see for instance \cite[p.~43]{ChaZak97}):


\begin{Theorem}
Let $0 < n \in \omega$, and let $\A$ be a Heyting algebra.
\benroman
\item $\A$ has depth at most $n$ if and only if it satisfies the equation $d_{n} \thickapprox 1$,  where
\begin{align*}
d_{1} &\coloneqq x_{1} \lor ( x_{1} \to 0)\\
d_{n+1} &\coloneqq x_{n+1} \lor ( x_{n+1} \to d_{n}), \text{ for all }n  \geq 1.
\end{align*}
\item $\A$ has width at most $n$ if and only if it satisfies the equation $w_{n} \thickapprox 1$, where
\[
w_{n} \coloneqq \bigvee_{i = 0}^{n}\big( x_{i} \to \bigvee_{j \ne i}x_{j}\big).
\]
\eroman
As a consequence, both $\class{D}_{n}$ and $\class{W}_{n}$ are varieties.
\end{Theorem}

Note that (ii) generalizes the fact that a Heyting algebra is a G\"{o}del algebra if and only if it satisfies $(x \to y)\lor(y \to x) \approx 1$.
Notice also that finite algebras have bounded depth, so the following general result from \cite[Thm.\ 5.3]{BMR17} can be viewed as a generalization of Proposition~\ref{Prp:finitely-generated}.\footnote{Theorem \ref{Thm:ES-bounded-depth} has a strengthening in the setting of substructural logics \cite{MRW4}. 
For further information on variants of the Beth property in substructural logics, the reader may consult \cite{Hoo00,KO10,Mo06a,Urq99}.}

\begin{Theorem}\label{Thm:ES-bounded-depth}
Let $0 < n \in \omega$. Every variety of Heyting algebras, whose members have depth at most $n$, has surjective epimorphisms.
\end{Theorem}
Kuznetsov showed that there is a continuum of varieties of Heyting algebras all of whose members have depth at most $3$ \cite{Kuz75,Kuznet74proc}. So, the result above supplies a continuum of varieties \emph{with} the ES property.

On the other hand, until now, only one \emph{ad hoc} example of a variety of Heyting algebras \emph{without} the ES property has been exhibited \cite[Cor.\ 6.2]{BMR17}. This variety is generated by an algebra that we call $(\D_2^{\infty})^{\ast}$, which has width at most $2$. This prompted the question of whether the variety $\class{W}_{2}$ (or $\class{W}_{n}$ in general) has the ES property or not. We shall settle this question after introducing a ``summing'' construction, which will be used to build $(\D_2^{\infty})^{\ast}$ as well as other algebras that have proper epic subalgebras.

\section{Sums of Heyting algebras}
\label{sec:sums}

Let $\A$ and $\B$ be Heyting algebras. The \emph{sum} $\A + \B$ 
is the Heyting algebra obtained by pasting $\B$ \emph{below} $\A$, gluing the top element of $\B$ to the bottom element of $\A$. To give a more formal definition, it is convenient to assume that the universes of $\A$ and $\B$ are disjoint. Moreover, let us denote by $\leq^{\A}$ and $\leq^{\B}$ the lattice orders of $\A$ and $\B$ respectively. Then $\A + \B$ is the unique Heyting algebra with universe $(A \smallsetminus \{ 0^{\A} \}) \cup B$ whose lattice order $\leq$ is defined as follows: for every $a, b \in (A \smallsetminus \{ 0^{\A} \}) \cup B$,
\begin{align*}
b \leq a \Longleftrightarrow& \text{ either }(a, b \in A \text{ and }b \leq^{\A}a) \text{ or }(a, b \in B \text{ and }b \leq^{\B}a)\\
& \text{ or }(b \in B \text{ and }a \in A).
\end{align*}
As $+$ is clearly associative, there is no ambiguity in writing $\A_{1} + \dots + \A_{n}$ for the descending chain of finitely many Heyting algebras $\A_{1}, \dots, \A_{n}$, each glued to the previous one.

To obtain interesting results 
about epimorphisms in varieties of Heyting algebras that are not consequences of
Theorem~\ref{Thm:ES-bounded-depth}, we will need to consider Heyting algebras with unbounded depth. It is therefore useful to introduce an infinite generalization of this construction.
Let $\{ \A_{n} \colon n \in \omega \}$ be a family of Heyting algebras with disjoint universes, and let $\bot$ be a fresh element. The \emph{sum} $\sum \A_{n}$ is the unique Heyting algebra with universe
\[
\{ \bot \} \cup \bigcup_{n \in \omega}(A_{n} \smallsetminus \{ 0^{\A_{n}} \})
\]
and whose lattice order is defined as follows: for every $a, b \in \sum \A_{n}$,
\begin{align*}
a \leq b \Longleftrightarrow& \text{ either }a = \bot \text{ or } (a, b \in A_{n} \text{ for some }n \in \omega\text{ and }a \leq^{\A_{n}}b)\\
& \text{ or }(a \in A_{n} \text{ and }b \in A_{m} \text{ for some }n, m \in \omega \text{ such that }n > m).
\end{align*}
In words, $\sum \A_n$ is a tower of algebras, each pasted below the previous, with a new bottom element. When $\{ \A_{n} \colon n \in \omega \}$ is a family consisting of copies of the same algebra $\A$, we write $\A^{\infty}$ instead of $\sum \A_{n}$.

Sums of Heyting algebras have found various applications in the study of intermediate logics. 
See for instance \cite{BBdeJ08,GeKuz70,KuzGer70a,Mak77a}, but note that in the usual definition, subsequent algebras are added on top, instead of below. For finitely many summands this difference is immaterial.

For present purposes, it is convenient to describe the dual spaces of sums of Heyting algebras as well. Let $\mathbb{X} = \langle X; \leq^{\mathbb{X}}\rangle$ and $\mathbb{Y} = \langle Y; \leq^{\mathbb{Y}}\rangle$ be two posets (with disjoint universes). 
Their \emph{sum} $\mathbb{X} + \mathbb{Y}$ is the poset with universe $X \cup Y$ and whose order relation $\leq$ is defined as follows: for every $x, y \in X \cup Y$,
\begin{align*}
x \leq y \Longleftrightarrow& \text{ either }(x, y \in X \text{ and }x \leq^{\mathbb{X}}y) \text{ or }(x, y \in Y \text{ and }x \leq^{\mathbb{Y}}y)\\
& \text{ or }(x \in X \text{ and }y \in Y).
\end{align*}
So, $\mathbb{X} + \mathbb{Y}$ is the poset obtained by placing $\mathbb{Y}$ \emph{above} $\mathbb{X}$.
Then let $\{ \mathbb{X}_{n} \colon n \in \omega \}$ be a family of posets with disjoint universes, and let $\top$ be a fresh element. The \emph{sum} $\sum \mathbb{X}_{n}$ is the poset with universe
\[
\{ \top \} \cup \bigcup_{n \in \omega} X_{n}
\]
and order relation $\leq$ defined as follows: for every $x, y \in \sum X_{n}$,
\begin{align*}
x \leq y \Longleftrightarrow& \text{ either }y = \top \text{ or } (x, y \in X_{n} \text{ for some }n \in \omega\text{ and }x \leq^{\mathbb{X}_{n}}y)\\
& \text{ or }(x \in X_{n} \text{ and }y \in X_{m} \text{ for some }n, m \in \omega \text{ such that }n < m).
\end{align*}
Hence, $\sum \mathbb{X}_n$ is obtained by placing each successive poset above the previous and adding a new top element.

Now, let $\boldsymbol{X}$ and $\boldsymbol{Y}$ be two Esakia spaces with disjoint universes. 
The \emph{sum} $\boldsymbol{X} + \boldsymbol{Y}$ is the Esakia space, whose underlying poset is $\langle X; \leq^{\boldsymbol{X}}\rangle + \langle Y; \leq^{\boldsymbol{Y}}\rangle$, endowed with the topology consisting of the sets $U \subseteq X \cup Y$ such that $U \cap X$ and $U \cap Y$ are open, respectively, 
in $\boldsymbol{X}$ and $\boldsymbol{Y}$.

Similarly, let $\{ \boldsymbol{X}_{n} \colon n \in \omega \}$ be a family of Esakia spaces with disjoint universes. The \emph{sum} $\sum \boldsymbol{X}_{n}$ is the Esakia space, whose underlying poset is $\sum\langle X_{n}; \leq^{\boldsymbol{X}_{n}}\rangle$, equipped with the topology
\begin{align*}
\tau = \?\? & \{ U \colon U \cap X_{n} \text{ is open in }\boldsymbol{X}_{n} \text{ for all }n \in \omega, \text{ and }\\
& \text{ if }\top \in U \text{, then there exists }n \in \omega \text{ with }\bigcup_{n \leq m} X_{m} \subseteq U \}.
\end{align*}
When $\{ \boldsymbol{X}_{n} \colon n \in \omega \}$ consists of copies of the same Esakia space $\boldsymbol{X}$, we write $\boldsymbol{X}^{\infty}$ instead of $\sum \boldsymbol{X}_{n}$.

\begin{Lemma}
\label{Lem:SumDuality}
Let $\{ \A, \B \} \cup \{ \A_{n} \colon n \in \omega \}$ be a family of Heyting algebras. The Esakia spaces $(\A + \B)_{\ast}$ and $(\sum\A_{n})_{\ast}$ are isomorphic, respectively, to $\A_{\ast} + \B_{\ast}$ and $\sum\A_{n \ast}$.
\end{Lemma}

\begin{proof}
We sketch the proof only for the case of $\sum\A_{n}$. Observe that the set of all non-zero elements of $\sum\A_{n}$ forms a prime filter $G_{0}$. Keeping this in mind, we define a map
\[
f \colon \sum\A_{n \ast} \to (\sum\A_{n})_{\ast}\,, 
\]
setting $f(\top) \coloneqq G_{0}$ and for every $n \in \omega$ and $F \in \Pr \A_{n}$,
\[
f(F) \coloneqq \{ a \in \sum\A_{n} \colon a \geq b \text{ for some } b \in F \}.
\]
It is not difficult to see that $f$ is bijective and order-preserving.  Therefore it only remains to prove that $f$ is continuous and satisfies (\ref{Eq:Esakia-morphism}).

To show the latter, suppose that $\top \neq F \in \sum\A_{n \ast}$ and $G_{0} \neq G \in (\sum\A_{n})_{\ast}$ such that $f(F) \leq G$. 
Then $F \in \Pr \A_{j}$ for some $j \in \omega$. We may therefore let $k \in \omega$ be the least $k \geq j$ such that $G \cap A_{m} = \emptyset$ for every $m \geq k$. If we let $G' = G \cap A_{k}$, it follows that $G' \in \Pr \A_{k} \subseteq \sum\A_{n \ast}$, with $F \leq G'$ and $f(G') = G$.

To prove that $f$ is continuous, first consider some subbasic clopen of the form $\gamma^{\sum\A_{n}}(a)$ with $a \in \sum\A_{n} \smallsetminus \{ 0 \}$, see (\ref{Eq:gamma}) if necessary. We have 
$a \in A_{k}$ for some $k \in \omega$ and, therefore,
\[
f^{-1}[\gamma^{\sum\A_{n}}(a)] = \{ \top \} \cup \gamma^{\A_{k}}(a) \cup \bigcup_{m > k} A_{m \ast}.
\]
Clearly the sets $\gamma^{\A_{k}}(a)$ and $\{ \top \} \cup \bigcup_{m > k} A_{m \ast}$ are open in $\sum\A_{n \ast}$. Now, similarly, consider a subbasic clopen of the form $\gamma^{\sum\A_{n}}(a)^{c}$. Then
\[
f^{-1}[\gamma^{\sum\A_{n}}(a)^{c}] = f^{-1}[\gamma^{\sum\A_{n}}(a)]^{c} = \gamma^{\A_{k}}(a)^{c} \cup \bigcup_{m = 0}^{k-1} A_{m \ast},
\]
which is also clearly open in $\sum\A_{n \ast}$. This shows that $f$ is continuous.
Therefore, $f$ is a bijective Esakia morphism, and hence an isomorphism.
\end{proof}

As an example, we can now construct $(\D_{2}^{\infty})^{\ast}$, which witnesses the failure of the ES property in the variety generated by it. 
We let $\D_2$ be the Esakia space with two incomparable elements. It follows that ${\D_2}^{\ast} \cong \2 \times \2$, where $\2$ is the 2-element Boolean algebra. From Lemma~\ref{Lem:SumDuality}, it follows that $(\D_{2}^{\infty})^{\ast} \cong (\2 \times \2)^{\infty}$, i.e., it is $\omega$ copies of the 4-element diamond, each pasted below the previous one, with a new bottom element.%
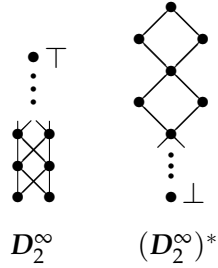
\footnote{The algebra used in \cite[Cor.\ 6.2]{BMR17} to generate a variety without the ES property is in fact $\2 + (\D_{2}^{\infty})^{\ast}$, which has the virtue of being FSI. For us, however, this difference is immaterial, since $\2 + (\D_{2}^{\infty})^{\ast}$ and $(\D_{2}^{\infty})^{\ast}$ generate the same variety, and both have epic subalgebras in this variety.}

\begin{center}
\begin{tabular}{ccc}

\begin{picture}(17,61)(17,0)
%
%

\put(30,52){$\top$}
\put(26,56){\circle*{4}}

\put(26,49){\circle*{2}}
\put(26,44){\circle*{2}}
\put(26,39){\circle*{2}}

\put(20,27){\circle*{4}}
\put(20,27){\line(0,1){5}}
\qbezier(20,27)(20,27)(25,32)

\put(32,27){\circle*{4}}
\put(32,27){\line(0,1){5}}
\qbezier(32,27)(32,27)(27,32)

\put(32,15){\circle*{4}}
\put(32,27){\line(0,-1){12}}
\put(20,27){\line(1,-1){12}}

\put(20,15){\circle*{4}}
\put(20,15){\line(0,1){12}}
\put(20,15){\line(1,1){12}}

\put(20,3){\circle*{4}}
\put(20,3){\line(0,1){12}}
\put(20,3){\line(1,1){12}}

\put(32,3){\circle*{4}}
\put(32,15){\line(0,-1){12}}
\put(20,15){\line(1,-1){12}}
\end{picture}

&&

\begin{picture}(50,93)(12,25)
%
%

\put(40,26){$\bot$}
\put(36,28){\circle*{4}}

\put(36,36){\circle*{2}}
\put(36,40){\circle*{2}}
\put(36,44){\circle*{2}}

\qbezier(36,52)(36,52)(41,47)
\qbezier(36,52)(36,52)(31,47)
\put(36,52){\circle*{4}}
\put(36,52){\line(1,1){12}}
\put(36,52){\line(-1,1){12}}
\put(48,64){\circle*{4}}
\put(24,64){\circle*{4}}
\put(36,76){\line(1,-1){12}}
\put(36,76){\line(-1,-1){12}}
\put(36,76){\circle*{4}}

\put(36,76){\line(1,1){12}}
\put(36,76){\line(-1,1){12}}
\put(48,88){\circle*{4}}
\put(24,88){\circle*{4}}
\put(36,100){\line(1,-1){12}}
\put(36,100){\line(-1,-1){12}}
\put(36,100){\circle*{4}}

\end{picture}
\vspace{0.5em}
\\
$\D_{2}^{\infty}$ &\hspace{5.0em}& $(\D_{2}^{\infty})^{\ast}$
\end{tabular}
\end{center}

\section{Varieties of bounded width}

As we mentioned, the variety of G\"odel algebras $\class{W}_{1}$ is known to have the ES property. On the other hand, in this section we show that for every $n \geq 2$, the variety $\class{W}_{n}$ lacks the ES property. To this end, we will rely on the following technical observation:

\begin{Lemma}\label{Lem:trick-width}
Let $0 < n \in \omega$ and let $f \colon \boldsymbol{Y} \to \boldsymbol{X}$ be an Esakia morphism between Esakia spaces of width at most $n$ such that
\benroman
\item $\boldsymbol{Y}$ has a minimum $\bot$, and 
\item for every $z \in X$ different from the maximum of $\boldsymbol{X}$ (if any), if $f(\bot) < z$, then there is an antichain of $n$ elements in ${\uparrow} f(\bot)$, which contains $z$.
\eroman
Also, let ${\uparrow} f(\bot)^{\top}$ be the upset ${\uparrow} f(\bot)$ without the maximum of $\boldsymbol{X}$ (if any). Then there is a subposet $\langle Z; \leq^{\boldsymbol{Y}} \rangle$ of $\boldsymbol{Y}$ such that the restriction
\[
f \colon \langle Z; \leq^{\boldsymbol{Y}} \rangle \to \?\? \langle{\uparrow} f(\bot)^{\top}; \leq^{\boldsymbol{X}}\rangle
\]
is a poset isomorphism.
\end{Lemma}

\begin{proof}
Observe that, since $f$ is an Esakia morphism, ${\uparrow} f(\bot)^{\top}$ coincides with $f[Y]$ without the maximum of $\bs{X}$ (if any). 
Suppose
$z \in  {\uparrow} f(\bot)^{\top} \smallsetminus \{ f(\bot) \}$ and define
\[
T_{z} \coloneqq \{ a \in Y \colon f(a) = z \}.
\]
Observe that $T_{z} \ne \emptyset$.

We claim that $T_{z}$ is a chain in $\boldsymbol{X}$. 
Indeed, since $f(\bot) <^{\boldsymbol{X}} z$,
assumption (ii) shows that 
$z$ belongs to an antichain $\{ x_{1}, \dots, x_{n-1}, z \}$ of $n$ elements in ${\uparrow} f(\bot)$. Since $f$ is an Esakia morphism, there are $y_{1}, \dots, y_{n-1} \in Y$ such that $f(y_{i}) = x_{i}$ for $i = 1, \dots, n-1$. Together with the fact that $f$ is order-preserving, this implies that the set $\{ y_{1}, \dots, y_{n-1}, a \}$ is an antichain of $n$ elements in $\boldsymbol{Y}$, for every $a \in T_{z}$. 

Now, suppose with a view to obtaining a contradiction 
that $T_{z}$ is not a chain. Then there are two incomparable elements $a, c \in T_{z}$. Hence $\{ y_{1}, \dots, y_{n-1}, a, c \}$ is an antichain of $n+1$ elements in $\boldsymbol{Y}$. Together with the fact that $\boldsymbol{Y}$ has a minimum element by assumption (i), we conclude that $\boldsymbol{Y}$ does not have width at most $n$. But this contradicts the assumptions, thus establishing the claim.

By Lemma \ref{Lem:correspondences}(\ref{Lem:correspondences:Chains}), the chain $T_{z}$ has a maximum element, which we denote by $\max(T_{z})$. Consider the set
\[
Z \coloneqq \{ \max(T_{z}) \colon z \in {\uparrow} f(\bot)^{\top} \smallsetminus \{ f(\bot) \} \} \cup \{ \bot \}.
\]
Clearly $Z \subseteq Y$, and it is easy to verify that the restriction
\[
f \colon \langle Z; \leq^{\boldsymbol{Y}} \rangle \to \langle {\uparrow} f(\bot)^{\top}; \leq^{\boldsymbol{X}}\rangle
\]
is a surjective order-preserving map. In order to prove that $f$ is a poset isomorphism, it remains only to show that $f$ is order-reflecting. To this end, consider $z_{1}, z_{2} \in Z$ such that $f(z_{1}) \leq^{\boldsymbol{X}} f(z_{2})$. If $z_{1} = \bot$, then clearly $z_{1} = \bot \leq^{\boldsymbol{Y}} z_{2}$, and we are done. So, consider the case where $z_{1} \ne \bot$.  In particular, this implies that $f(z_{1}),f(z_{2}) \in {\uparrow} f(\bot)^{\top} \smallsetminus \{ f(\bot) \}$, and therefore that $z_{2} = \max(T_{f(z_{2})})$. Now, since $f(z_{1}) \leq^{\boldsymbol{X}} f(z_{2})$ and $f$ is an Esakia morphism, there exists $z_{3} \in Y$ such that $z_{1} \leq^{\boldsymbol{Y}} z_{3}$ and $f(z_{3}) = f(z_{2})$. But then $z_{3} \in T_{f(z_{2})}$, which implies that $z_{3} \leq^{\boldsymbol{Y}} z_{2}$ and, therefore, that $z_{1} \leq^{\boldsymbol{Y}} z_{3} \leq^{\boldsymbol{Y}} z_{2}$. Thus, we conclude that $f$ is order-reflecting, as desired.
\end{proof}

For $1 < n \in \omega$, let $\mathbb{X}_{n} = \langle X_{n}; \leq \rangle$ be the poset whose universe 
is $\{ a_{1}, \dots, a_{n}, b_{1}, \dots, b_{n} \}$, and whose order relation is defined as follows (see the picture below): for 
all 
$x, y \in X_{n}$,
\begin{align*}
x \leq y \Longleftrightarrow &\text{ either }x = y \text{ or }(x= a_{1} \text{ and }y \in \{ b_{2}, \dots, b_{n} \})\\
& \text{ or }(x = a_{m} \text{ for some }m > 1\text{ and }y \in \{ b_{1}, b_{m} \}).
\end{align*}

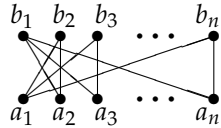
\begin{center}

\begin{picture}(80,47)
%
%


\put(5,10){\usebox{\XnPic}}

{
\small
\put(0,2){$a_1$}
\put(0,39){$b_1$}
\put(16,2){$a_2$}
\put(16,39){$b_2$}
\put(32,2){$a_3$}
\put(32,39){$b_3$}
\put(70,2){$a_n$}
\put(70,39){$b_n$}
}
\end{picture}

\end{center}

Also, let $\boldsymbol{X}_{n}$ be the Esakia space obtained by endowing $\mathbb{X}_{n}$ with the discrete topology. 

\begin{Remark}
\label{Rem:duals-of-Xn}
The reader may wonder whether there is an intelligible way to understand the algebraic duals of the spaces $\boldsymbol{X}_{n}$. Although we will not rely on this observation, it is possible to show that the lattice reduct $\C$ of ${\boldsymbol{X}_{n}}^{\ast}$ is the distributive lattice obtained as follows. Let $\A$ be the Boolean lattice of $2^{n-1}$ elements, extended with two new extrema. Moreover, let $\B$ be the three-element chain. Then $\C$ is obtained by computing the coproduct of $\A$ and $\B$ in the category of bounded distributive lattices and, subsequently, removing the two extrema from it. 

An easy proof of this fact can be given by means of Priestley duality \cite{Pr70,Pr72} for bounded distributive lattices, in which finite coproducts of finite algebras correspond to direct products of posets in the natural sense \cite[Thm.\ 3, Sec.\ VII.7]{BaDw74}.
\qed
\end{Remark}

The next observation follows immediately from the definitions:

\begin{Lemma}\label{Lem:easy-bounded-width}
For $1 < n \in \omega$, the Esakia space $\boldsymbol{X}_{n}^{\infty}$ has width at most $n$.
\end{Lemma}

Now, let us introduce a notation for referring to the elements of $\boldsymbol{X}_{n}^{\infty}$. Recall that the universe of $\boldsymbol{X}_{n}^{\infty}$ is a chain of copies of $X_{n}$ of order type $\omega$, plus an new maximum element $\{ \top \}$. We denote the elements of the lowest component of $X_{n}$ in $\boldsymbol{X}_{n}^{\infty}$ as follows:
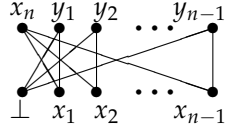
\begin{center}

\begin{picture}(83,47)
%
%


\put(5,12){\usebox{\XnPic}}

{
\small
\put(0,2){$\bot$}
\put(0,41){$x_n$}
\put(16,2){$x_1$}
\put(16,41){$y_1$}
\put(32,2){$x_2$}
\put(32,41){$y_2$}
\put(62,2){$x_{n-1}$}
\put(62,41){$y_{n-1}$}
}
\end{picture}

\end{center}
Similarly, for $1 \leq k \in \omega$, we denote the elements of the $(k+1)$-th component of $X_{n}$ in $\boldsymbol{X}_{n}^{\infty}$ as follows:
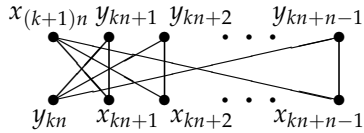
\begin{center}

\begin{picture}(125,47)(-2,0)
%
%


\put(10,10){\scalebox{1.5}[1]{\usebox{\XnLines}}}

\put(10,10){\circle*{4}}
\put(10,34){\circle*{4}}
\put(31,34){\circle*{4}}
\put(31,10){\circle*{4}}

\put(52,34){\circle*{4}}
\put(52,10){\circle*{4}}

\put(118,34){\circle*{4}}
\put(118,10){\circle*{4}}

\put(75,10){\circle*{2}}
\put(84,10){\circle*{2}}
\put(93,10){\circle*{2}}
\put(75,34){\circle*{2}}
\put(84,34){\circle*{2}}
\put(93,34){\circle*{2}}

{
\small
\put(2,2){$y_{kn}$}
\put(-7,41){$x_{(k+1)n}$}
\put(26,2){$x_{kn+1}$}
\put(26,41){$y_{kn+1}$}
\put(54,2){$x_{kn+2}$}
\put(54,41){$y_{kn+2}$}
\put(93,2){$x_{kn+n-1}$}
\put(93,41){$y_{kn+n-1}$}
}
\end{picture}

\end{center}
Keeping this in mind, let $R_{n}$ be the equivalence relation on $X_{n}^{\infty}$ defined as follows: for every $a, b \in X_{n}^{\infty}$,
\begin{align*}
\langle a, b \rangle \in R_{n} \Longleftrightarrow& \text{ either } a = b \text{ or}\\
&\text{ there is $k \in \omega$ such that $\{ a, b \} = \{ x_{k}, y_{k}\}$.}
\end{align*}
A pictorial rendering of the partition corresponding to 
the relation $R_{n}$ is given below:
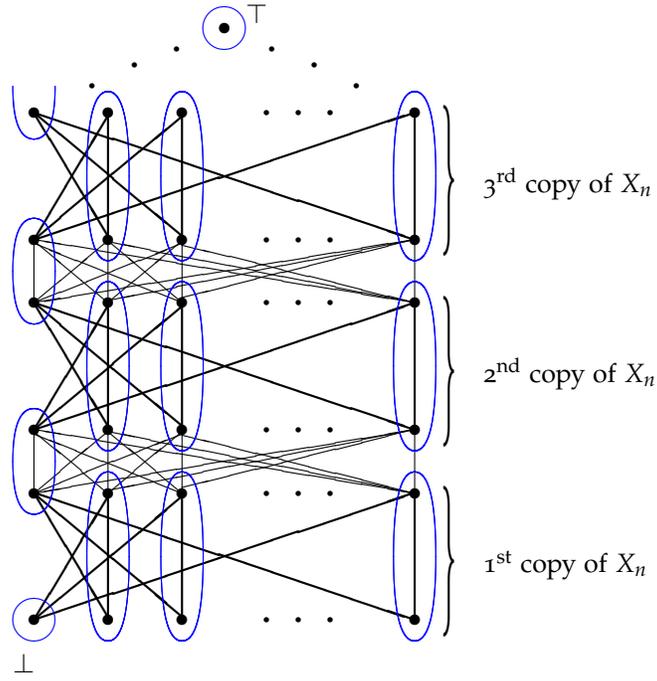
\begin{center}
\setlength{\unitlength}{2pt}
\begin{picture}(125,128)(0,7)
%
%


\newsavebox{\XnPicBig}
\savebox{\XnPicBig}(72,24)[bl]{
\thicklines
\put(0,0){\line(3,5){14}}
\put(0,24){\line(3,-5){14}}
\put(14,0){\line(0,1){24}}

\put(0,0){\line(6,5){28}}
\put(0,24){\line(6,-5){28}}
\put(28,0){\line(0,1){24}}

\put(0,0){\line(3,1){72}}
\put(0,24){\line(3,-1){72}}
\put(72,0){\line(0,1){24}}

\thinlines

\put(0,0){\circle*{2}}
\put(0,24){\circle*{2}}
\put(14,24){\circle*{2}}
\put(14,0){\circle*{2}}

\put(28,24){\circle*{2}}
\put(28,0){\circle*{2}}

\put(72,24){\circle*{2}}
\put(72,0){\circle*{2}}

\put(44,0){\circle*{1}}
\put(50,0){\circle*{1}}
\put(56,0){\circle*{1}}
\put(44,24){\circle*{1}}
\put(50,24){\circle*{1}}
\put(56,24){\circle*{1}}
}

\put(5,18){\usebox{\XnPicBig}}
\put(5,54){\usebox{\XnPicBig}}
\put(5,90){\usebox{\XnPicBig}}

\newsavebox{\ConnectLines}
\savebox{\ConnectLines}(72,12)[bl]{%
\put(0,0){\line(0,1){12}}
\put(0,0){\line(6,5){14}}
\put(0,0){\line(5,2){28}}

\put(14,0){\line(0,1){12}}
\put(14,0){\line(6,5){14}}
\put(14,0){\line(-6,5){14}}

\put(28,0){\line(0,1){12}}
\put(28,0){\line(-6,5){14}}
\put(28,0){\line(-5,2){28}}

\put(72,0){\line(0,1){12}}
\put(72,0){\line(-4,1){44}}
\put(72,0){\line(-5,1){58}}
\put(72,0){\line(-6,1){72}}

\put(72,12){\line(-4,-1){44}}
\put(72,12){\line(-5,-1){58}}
\put(72,12){\line(-6,-1){72}}
}
\put(5,42){\usebox{\ConnectLines}}
\put(5,78){\usebox{\ConnectLines}}

\put(50,126){\circle*{1}}
\put(58,123){\circle*{1}}
\put(66,119){\circle*{1}}
\put(32,126){\circle*{1}}
\put(24,123){\circle*{1}}
\put(16,119){\circle*{1}}

\put(41,130){\circle*{2}}


{\color{blue}

\put(41,130){\circle{8}}
\put(5,18){\circle{8}}

\newsavebox{\EqLong}
\savebox{\EqLong}(8,32){%
\qbezier(4,0)(0,0)(0,16)
\qbezier(0,16)(0,32)(4,32)
\qbezier(4,0)(8,0)(8,16)
\qbezier(8,16)(8,32)(4,32)
}

\put(15,14){\usebox{\EqLong}}
\put(29,14){\usebox{\EqLong}}
\put(73,14){\usebox{\EqLong}}
\put(15,50){\usebox{\EqLong}}
\put(29,50){\usebox{\EqLong}}
\put(73,50){\usebox{\EqLong}}
\put(15,86){\usebox{\EqLong}}
\put(29,86){\usebox{\EqLong}}
\put(73,86){\usebox{\EqLong}}

\newsavebox{\EqShort}
\savebox{\EqShort}(8,20){%
\qbezier(4,0)(0,0)(0,10)
\qbezier(0,10)(0,20)(4,20)
\qbezier(4,0)(8,0)(8,10)
\qbezier(8,10)(8,20)(4,20)
}
\put(1,38){\usebox{\EqShort}}
\put(1,74){\usebox{\EqShort}}

\newsavebox{\EqHalf}
\savebox{\EqHalf}(8,20){%

\qbezier(4,0)(0,0)(0,10)
\qbezier(4,0)(8,0)(8,10)
}
\put(1,104){\usebox{\EqHalf}}

}

\put(82,22){\scalebox{1}[5]{$\}$}}
\put(82,58){\scalebox{1}[5]{$\}$}}
\put(82,94){\scalebox{1}[5]{$\}$}}

{
\small
\put(1,8){$\bot$}
\put(90,27){1\textsuperscript{st} copy of $X_n$}
\put(90,63){2\textsuperscript{nd} copy of $X_n$}
\put(90,99){3\textsuperscript{rd} copy of $X_n$}
\put(45,131){$\top$}
}
\end{picture}
\setlength{\unitlength}{1pt}
\end{center}

\begin{Lemma}\label{Lem:correct-partition}
$R_{n}$ is a correct partition on $\boldsymbol{X}_{n}^{\infty}$, for $1 < n \in \omega$.
\end{Lemma}

\begin{proof}
The fact that $R_{n}$ is an equivalence relation satisfying condition (i) in the definition of correct partitions can be verified by inspecting the figure above. 

We turn now to condition (ii). First we claim that every finite set $U \subseteq X_{n}^{\infty} \smallsetminus \{ \top \}$ is clopen. The fact that $U$ is closed follows from the fact that, in Hausdorff spaces, every singleton is closed. On the other hand, the fact that $U$ is open follows from the definition of $\boldsymbol{X}_{n}^{\infty}$, together with the fact that $\boldsymbol{X}_{n}$ is endowed with the discrete topology. This establishes the claim.

Now, consider two distinct points $x, y \in X_{n}^{\infty}$ such that $\langle x, y \rangle \notin R_{n}$. We need to find a clopen $U$ such that $x \in U$ and $y \notin U$, which moreover is a union of equivalence classes of $R_{n}$. If $x \ne \top$, the equivalence class $x/ R_{n}$ is a finite subset of $X_{n}^{\infty} \smallsetminus \{ \top \}$. From the claim it follows that $x/ R_{n}$ is clopen. Since $y \notin x/ R_{n}$, we let $U \coloneqq x / R_{n}$ and we are done. Then consider the case where $x = \top$. We know that $y$ appears in the $k$-th component of $X_{n}$ in $\boldsymbol{X}_{n}^{\infty}$, for some $k \in \omega$. Consider the set $U \coloneqq {\uparrow} x_{(k+1)n}$. Clearly $x  = \top \in U$ and $y \notin U$. Moreover, $U$ is a union of equivalence classes of $R_{n}$. Finally, from the claim it follows that $U^{c}$ is clopen and, therefore, so is $U$.
\end{proof}

We are now ready to prove the main result of this section:

\begin{Theorem}\label{Thm:Es-bounded-width}
Let $1 < n \in \omega$, and let $\mathsf{K}$ be a variety of Heyting algebras, whose members have width at most $n$. If $\boldsymbol{X}_{n}^{\infty} \in \mathsf{K}_{\ast}$, then $\mathsf{K}$ lacks the ES property.
\end{Theorem}

\begin{proof}

Fix $1 < n \in \omega$, and let $\class{K}$ be a variety of Heyting algebras whose members have width at most $n$, such that $\boldsymbol{X}_{n}^{\infty} \in \class{K}_{\ast}$. To show that the ES property fails in $\class{K}$, we will employ Lemma~\ref{Lem:epic-subalgebras-dual}. In particular, by Lemma~\ref{Lem:correct-partition}, we know that $R_{n}$ is a correct partition of $\boldsymbol{X}_{n}^{\infty}$ that is clearly different from the identity relation. We suppose, with a view to obtaining a contradiction, 
that there exist $\bs{Y} \in \class{K}_{\ast}$ and a pair of \emph{different} homomorphisms $g,h \colon \boldsymbol{Y} \to \boldsymbol{X}_{n}^{\infty}$, such that $\langle f(y), g(y) \rangle \in R_{n}$ for every $y \in Y$.

We claim that there is an element $\bot \in Y$ and $0 < k \in \omega$ such that $\{ f(\bot), g(\bot) \} = \{ x_{kn}, y_{kn} \}$. First observe that, since $f \ne g$, there is $y \in Y$ such that $f(y) \ne g(y)$ and $\langle f(y), g(y)\rangle \in R_{n}$. Then the equivalence class $f(y)/R_{n}$ is not a singleton. This implies that $\{ f(y), g(y) \} = \{ x_{m}, y_{m} \}$ for some positive integer $m$.

If $m$ is a multiple of $n$, we are done. Consider the case where $m = tn + s$ for some $s,t \in \omega$, where $0 < s < n$. Suppose without loss of generality that  $f(y) = x_{m}$. Then we have that
\[
f(y) = x_{m} \leq^{\boldsymbol{X}^{\infty}_{n}} x_{(t+1)n}.
\]
Since $f$ is an Esakia morphism, there is an element $x_{(t+1)n}' \in Y$ such that
\[
y \leq^{\boldsymbol{Y}} x_{(t+1)n}' \text{ and }f(x_{(t+1)n}') = x_{(t+1)n}.
\]
Now, recall that $g(y) = y_{m}$. In particular, this implies that
\[
y_{m} \leq^{\boldsymbol{X}_{n}^{\omega}} g(y) \leq^{\boldsymbol{X}_{n}^{\infty}}  g(x_{(t+1)n}').
\]
Together with the fact that $y_{m} \nleq^{\boldsymbol{X}_{n}^{\infty}} x_{(t+1)n}$, this implies that $g(x_{(t+1)n}') \ne x_{(t+1)n}$. Now, from the fact that $\langle f(x_{(t+1)n}'), g(x_{(t+1)n}') \rangle \in R_{n}$ and $f(x_{(t+1)n}') = x_{(t+1)n}$, it follows that $g(x_{(t+1)n}') = y_{(t+1)n}$. Hence, setting $\bot \coloneqq x_{(t+1)n}'$ and $k \coloneqq t+1$, we have that $0 < k \in \omega$ and $\{ f(\bot), g(\bot) \} = \{ x_{kn}, y_{kn} \}$, establishing the claim.

Now, let $\bot \in Y$ be the element given by the claim. From the definition of an Esakia space we know that the upset ${\uparrow} \bot$ is closed and, therefore, an E-subspace of $\boldsymbol{Y}$. Moreover, ${\uparrow} \bot_{\ast} \in \class{K}_{\ast}$, since $\class{K}$ is a variety.  Finally, by Lemma \ref{Lem:correspondences}(\ref{Lem:correspondences:ImagesAreSubspaces}), the restrictions of $f$ and $g$ to ${\uparrow} \bot$ are Esakia morphisms. Therefore we can assume without of generality that $\boldsymbol{Y} =  {\uparrow} \bot$ (otherwise we replace $\boldsymbol{Y}$ with ${\uparrow} \bot$).

Recall from the claim that $\{ f(\bot), g(\bot) \} = \{ x_{kn}, y_{kn} \}$. We can assume without loss of generality that $f(\bot) = x_{kn}$ and that $g(\bot) = y_{kn}$. Then observe that the  Esakia spaces $\bs{Y}$ and $\boldsymbol{X}_{n}^{\infty}$, and the Esakia morphism $f \colon \boldsymbol{Y} \to \boldsymbol{X}_{n}^{\infty}$ satisfy the assumptions of Lemma \ref{Lem:trick-width}. Therefore, $\boldsymbol{Y}$ has a subposet $\langle Z; \leq^{\boldsymbol{Y}}\rangle$ such that the restriction
\[
f \colon \langle Z; \leq^{\boldsymbol{Y}}\rangle \to \langle {\uparrow} f(\bot)^{\top}; \leq^{\boldsymbol{X}_{n}^{\infty}}\rangle
\]
is a poset isomorphism. For the sake of simplicity, we denote the elements of $Z$ exactly as their alter egos in ${\uparrow} f(\bot)^{\top}$. 

Under this convention, we have that
\[
Z = \{ x_{kn + m} \colon m \in \omega \} \cup \{ y_{kn + m} \colon m \in \omega \}
\]
and that $f(x_{i}) = x_{i}$ and $f(y_{i}) = y_{i}$, for every $x_{i}, y_{i} \in Z$. 
(Note that $\bot$ has the label $x_{kn}$.)
On the other hand, since $g$ is order-preserving and $g(\bot) = y_{kn}$, we have $g(x_{i}) = g(y_{i}) = y_{i}$, for all $x_{i}, y_{i} \in Z$. Consequently, for all $x_{i}, y_{i} \in Z$,
\[
f(x_{i}) = x_{i} \text{ and }f(y_{i}) = g(x_{i}) = g(y_{i}) = y_{i}.
\]

Now, consider the set
\[
C \coloneqq \{ y_{(k+1)n}, x_{(k+1)n + 1}, x_{(k+1)n + 2}, \dots, x_{(k+1)n + n -1} \} \subseteq {\uparrow} f(\bot)^{\top}.
\]
Observe that $C$ is an antichain of $n$ elements in $\boldsymbol{X}_{n}^{\infty}$. Since $f$ is a poset isomorphism between $\langle Z; \leq^{\boldsymbol{Y}}\rangle$ and $\langle {\uparrow} f(\bot)^{\top}; \leq^{\boldsymbol{X}_{n}^{\infty}}\rangle$, we obtain that $C$ is also an antichain of $n$ elements in $\langle Z; \leq^{\boldsymbol{Y}}\rangle$ and, therefore, in $\boldsymbol{Y}$. Since $\boldsymbol{Y}$ has width at most $n$ and has a minimum element $\bot$, every element of $Y$ should be comparable with at least one member of $C$.

Now, observe that
\[
g(\bot) = y_{kn} \leq^{\boldsymbol{X}_{n}^{\infty}} x_{(k+2)n}.
\]
Since $g$ is an Esakia morphism, there exists $a \in Y$ such that $g(a) = x_{(k+2)n}$. Recall that for $1 \leq i \leq n-1$,
\[
g(x_{(k+1)n + i}) = y_{(k+1)n + i}.
\]
Since $y_{(k+1)n + i}$ is incomparable with $g(a) = x_{(k+2)n}$ in $\boldsymbol{X}_{n}^{\infty}$ and $g$ is order-preserving, we conclude that $a$ is incomparable with $x_{(k+1)n + i}$ in $\boldsymbol{Y}$. Hence $a$ is incomparable with every element of $C \smallsetminus \{ y_{(k+1)n} \}$ in $\boldsymbol{Y}$. On the other hand, we know that $a$ must be comparable with at least one element of $C$. As a consequence, $a$ must be comparable with $y_{(k+1)n}$ in $\boldsymbol{Y}$. Together with the fact that $g$ is order-preserving, this means that $g(a) = x_{(k+2)n}$ and $g(y_{(k+1)n}) = y_{(k+1)n}$ are comparable in $\boldsymbol{X}_{n}^{\infty}$, which is a contradiction.
\end{proof}

As a consequence of the above result, we obtain that for  $1 < n \in \omega$ the variety of all Heyting algebras of width at most $n$ lacks the ES property.

\begin{Corollary}\label{Cor:width-n}
For $1 < n \in \omega$, the members of the interval $[ \VVV(\boldsymbol{X}_{n}^{\infty\ast}), \class{W}_{n}]$ of the subvariety lattice of $\class{HA}$ lack the ES property. In particular, $\class{W}_{n}$  lacks the ES property.
\end{Corollary}

\begin{proof}
Immediate from Lemma \ref{Lem:easy-bounded-width} and Theorem \ref{Thm:Es-bounded-width}.
\end{proof}

At this stage the reader may wonder how many varieties of Heyting algebras lack the ES property. Later on (Theorem \ref{Thm:continuum}), we will show that there is already a continuum of locally finite varieties without the ES property whose members have width at most $2$.

\begin{Remark}
As we mentioned, the logical counterpart of the ES property is the infinite Beth property. One naturally asks, therefore, what failure of the infinite Beth property is related to the failure of the ES property in Theorem \ref{Thm:Es-bounded-width}. To clarify this point, observe that the proof of the theorem shows that the subalgebra of $(\boldsymbol{X}_{n}^{\infty})^{\ast}$ corresponding to the correct partition $R_{n}$ is proper and $\class{K}$-epic. In the cases where $n$ is either $2$ or $3$, the algebra $(\boldsymbol{X}_{n}^{\infty})^{\ast}$ is depicted below, where the encircled elements constitute the proper $\mathsf{K}$-epic subalgebra corresponding to $R_{n}$.

Let $\A \leq \B$ be 
one of the two epic situations depicted below. Observe that for $1 \leq k \in \omega$, there is a \emph{unique} element $a_{k} \in C_{k}$ which is incomparable with every element of $(A \cap C_{k}) \smallsetminus \{ 1 \}$ (see picture below). We call $a_{k}$ the \emph{sibling} of $(A \cap C_{k})$. Roughly speaking, Theorem \ref{Thm:Es-bounded-width} shows that siblings are \emph{implicitly}, but not \emph{explicitly}, definable in terms of $A$. Since $\B$ is generated by $A$ together with the siblings $\{ a_{k} \colon k \in \omega \}$, this observation implies that all the elements of $B$ are implicitly, but not explicitly, definable in terms of $A$, witnessing a failure of the infinite Beth property.\footnote{Formally speaking a failure of the infinite Beth property should be given in terms of two disjoint sets of variables $X$ and $Z$, and a set of formulas $\Gamma$ over $X \cup Z$. Our informal explanation can be amended, taking $X \coloneqq A$, $Z \coloneqq B \smallsetminus A$, and letting $\Gamma$ be the inverse image of $\{ 1^{\B}\}$ under the natural homomorphism from the term algebra $\Fm(X \cup Z)$ to $\B$.}

As a matter of fact, the above explanation is not restricted to the case where $n = 2, 3$. In general, let $\A$ be the subalgebra dual to $R_{n}$ on $\B \coloneqq (\boldsymbol{X}_{n}^{\infty})^{\ast}$. Whenever $1 \leq k \in \omega$, then $C_{k}$ has a unique element $a_{k}$ that is incomparable with all elements of $(A \cap C_{k}) \smallsetminus \{ 1 \}$. Indeed, $a_{k} = {\uparrow} y_{n(k-1)}$ when $k > 0$, and $a_{k} = {\uparrow} \bot$ otherwise.
Here, as before, the $a_{k}$'s are implicitly, but not explicitly, definable in terms of $A$ and, together with $A$, they generate the algebra $\B$, witnessing a failure of the infinite Beth property.
\qed
\end{Remark}

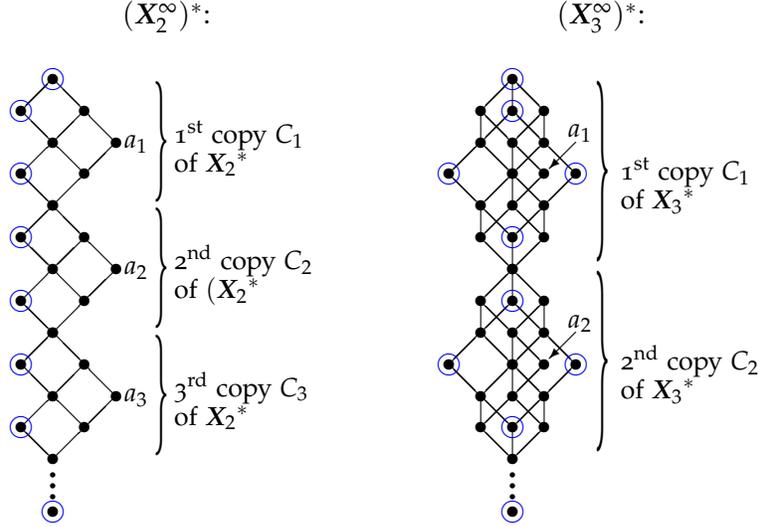
\begin{center}
\begin{tabular}{cp{2em}c}
$(\boldsymbol{X}_{2}^{\infty})^{\ast}$: && $(\boldsymbol{X}_{3}^{\infty})^{\ast}$:\vspace{1em} \\

\begin{picture}(120,175)
%
%


\newsavebox{\XtwoA}
\savebox{\XtwoA}(36,48){%

\put(0,12){\line(1,-1){12}}
\put(0,12){\circle*{4}}
\put(0,12){\line(1,1){24}}
\put(12,0){\circle*{4}}
\put(12,0){\line(1,1){24}}
\put(24,12){\circle*{4}}
\put(24,12){\line(-1,1){24}}
\put(36,24){\circle*{4}}
\put(36,24){\line(-1,1){24}}
\put(12,24){\circle*{4}}
\put(24,36){\circle*{4}}
\put(0,36){\line(1,1){12}}
\put(0,36){\circle*{4}}
\put(12,48){\circle*{4}}

{\color{blue}
\put(0,12){\circle{8}}
\put(0,36){\circle{8}}
}
}

\put(5,25){\usebox{\XtwoA}}
\put(5,73){\usebox{\XtwoA}}
\put(5,121){\usebox{\XtwoA}}

\put(17,5){\circle*{4}}

\put(17,11){\circle*{2}}
\put(17,15){\circle*{2}}
\put(17,19){\circle*{2}}

{\color{blue}
\put(17,5){\circle{8}}
\put(17,169){\circle{8}}
}

\put(55,38){\scalebox{1}[4]{$\}$}}
\put(55,86){\scalebox{1}[4]{$\}$}}
\put(55,134){\scalebox{1}[4]{$\}$}}

{
\small
\put(63,49){3\textsuperscript{rd} copy $C_3$}
\put(63,38){of ${\bs{X}_2}^{\ast}$}
\put(63,97){2\textsuperscript{nd} copy $C_2$}
\put(63,86){of $({\bs{X}_2}^{\ast}$}
\put(63,145){1\textsuperscript{st} copy $C_1$}
\put(63,134){of ${\bs{X}_2}^{\ast}$}
\put(44,47){$a_3$}
\put(44,95){$a_2$}
\put(44,143){$a_1$}
}
\end{picture}
&&
\begin{picture}(125,175)(0,72)
%
%


\newsavebox{\Cube}
\savebox{\Cube}(24,36){%

\put(12,0){\circle*{4}}
\put(12,0){\line(0,1){12}}
\put(12,12){\circle*{4}}
\put(12,0){\line(1,1){12}}
\put(24,12){\circle*{4}}
\put(12,0){\line(-1,1){12}}
\put(0,12){\circle*{4}}
\put(24,12){\line(0,1){12}}
\put(24,24){\circle*{4}}
\put(0,12){\line(0,1){12}}
\put(0,24){\circle*{4}}
\put(12,12){\line(1,1){12}}
\put(12,12){\line(-1,1){12}}
\put(24,12){\line(-1,1){12}}
\put(0,12){\line(1,1){12}}
\put(12,24){\circle*{4}}
\put(12,24){\line(0,1){12}}
\put(24,24){\line(-1,1){12}}
\put(0,24){\line(1,1){12}}
\put(12,36){\circle*{4}}
}

\newsavebox{\XtreeA}
\savebox{\XtreeA}(48,72){%

\put(12,0){\usebox{\Cube}}
\put(12,36){\usebox{\Cube}}

\put(0,36){\line(1,-1){12}}
\put(0,36){\circle*{4}}
\put(0,36){\line(1,1){12}}
\put(48,36){\line(-1,-1){12}}
\put(48,36){\circle*{4}}
\put(48,36){\line(-1,1){12}}
\put(36,36){\line(-1,-1){12}}
\put(36,36){\circle*{4}}
\put(36,36){\line(-1,1){12}}

{\color{blue}
\put(48,36){\circle{8}}
\put(0,36){\circle{8}}
\put(24,12){\circle{8}}
\put(24,60){\circle{8}}
}
}

\put(5,97){\usebox{\XtreeA}}
\put(5,169){\usebox{\XtreeA}}

\put(29,77){\circle*{4}}

\put(29,83){\circle*{2}}
\put(29,87){\circle*{2}}
\put(29,91){\circle*{2}}

{\color{blue}
\put(29,77){\circle{8}}
\put(29,241){\circle{8}}
}

\put(60,117){\scalebox{1}[6]{$\}$}}
\put(60,189){\scalebox{1}[6]{$\}$}}

{
\small
\put(70,131){2\textsuperscript{nd} copy $C_2$}
\put(70,120){of ${\bs{X}_3}^{\ast}$}
\put(70,203){1\textsuperscript{st} copy $C_1$}
\put(70,192){of ${\bs{X}_3}^{\ast}$}
\put(50,147){$a_2$}
\put(53,145){\vector(-1,-1){10}}
\put(50,219){$a_1$}
\put(53,217){\vector(-1,-1){10}}
}
\end{picture}
\end{tabular}
\end{center}

\section{The Kuznetsov-Ger\v{c}iu variety}

It is well-known that the free one-generated Heyting algebras is the \emph{Rieger-Nishimura} lattice $\bs{RN}$, depicted below \cite{Ni60,Ri49}. 
As a consequence, $\HHH(\bs{RN})$ is the class of one-generated Heyting algebras.

\begin{center}

\begin{picture}(62,144)
%
%
\put(38,3){\circle*{4}}
\put(38,3){\line(1,1){12}}
\put(38,3){\line(-1,1){12}}
\put(50,15){\circle*{4}}
\put(50,15){\line(-1,1){12}}
\put(26,15){\circle*{4}}
\put(26,15){\line(1,1){12}}
\put(38,27){\circle*{4}}
\put(12,5){$w_0$}
\put(50,6){$w_1$}
\put(42,25){$a_1$}

\put(26,15){\line(-1,1){12}}
\put(38,27){\line(-1,1){12}}
\put(14,27){\circle*{4}}
\put(14,27){\line(1,1){12}}
\put(26,39){\circle*{4}}
\put(0,17){$w_2$}
\put(30,37){$a_2$}

\put(38,27){\line(1,1){12}}
\put(50,39){\circle*{4}}
\put(50,39){\line(-1,1){12}}
\put(26,39){\line(1,1){12}}
\put(38,51){\circle*{4}}
\put(50,30){$w_3$}
\put(42,49){$a_3$}

\put(26,39){\line(-1,1){12}}
\put(38,51){\line(-1,1){12}}
\put(14,51){\circle*{4}}
\put(14,51){\line(1,1){12}}
\put(26,63){\circle*{4}}
\put(0,41){$w_4$}
\put(30,61){$a_4$}

\put(38,51){\line(1,1){12}}
\put(50,63){\circle*{4}}
\put(50,63){\line(-1,1){12}}
\put(26,63){\line(1,1){12}}
\put(38,75){\circle*{4}}
\put(50,54){$w_5$}
\put(42,73){$a_5$}

\put(26,63){\line(-1,1){12}}
\put(38,75){\line(-1,1){12}}
\put(14,75){\circle*{4}}
\put(14,75){\line(1,1){12}}
\put(26,87){\circle*{4}}
\put(0,65){$w_6$}
\put(30,85){$a_6$}

\put(38,75){\line(1,1){12}}
\put(50,87){\circle*{4}}
\put(50,87){\line(-1,1){10}}
\put(26,87){\line(1,1){10}}
\put(50,78){$w_7$}

\put(26,87){\line(-1,1){10}}

\put(32,105){\circle*{2}}
\put(32,111){\circle*{2}}
\put(32,117){\circle*{2}}

\put(32,129){\circle*{4}}
\put(36,127){$1$}
\end{picture}

%

\captionof{figure}{$\bs{RN}$}
\label{fig:RN}
\end{center}

The \emph{Kuznetsov-Ger\v{c}iu} variety is defined as follows:
\begin{equation}\label{Eq:generators-KG}
\class{KG} \coloneqq \VVV(\{ \A_{1} + \dots + \A_{n} \colon 0 < n \in \omega \text{ and } \A_{1}, \dots, \A_{n} \in \HHH(\bs{RN})\}).
\end{equation}
The variety $\class{KG}$ was introduced by Kuznetsov and Ger\v{c}iu \cite{GeKuz70,KuzGer70a} in their 
study of finite axiomatizability, and of the finite model property 
in varieties of Heyting algebras (also see \cite{BBdeJ08,BM19}).\ Remarkably, there is a continuum of subvarieties of $\class{KG}$ having the finite model property, and also a continuum lacking the finite model property \cite[Thm.\ 5.39(1), Cor.\ 5.41]{BBdeJ08}.
Notice that the variety $\VVV(\bs{RN})$ is contained in $\class{KG}$ and so are all of its subvarieties.

In the next sections we will provide a characterization of subvarieties of $\class{KG}$ with the ES property (Theorem \ref{Thm:ES-in-KG}). 

\begin{Lemma}[A.\ V.\ Kuznetsov and V.\ Ja.\ Ger\v{c}iu {\cite[Lem.\ 4]{KuzGer70a}}]\label{Lem:finite-FSI-KG}
If $\A \in \class{KG}$ is a finite FSI algebra, then $\A = \B_{1} + \dots + \B_{n}$ for some $\B_{1}, \dots, \B_{n} \in \HHH(\bs{RN})$.
\end{Lemma}



A variety is said to be \emph{locally finite} when its finitely generated members are finite. 
Recall that $\bs{2}$ denotes the two-element Boolean algebra.

\begin{Theorem}[N.\ Bezhanishvili, G.\ Bezhanishvili and D.\ de Jongh {\cite[Thms.\ 8.49 and 8.54]{BBdeJ08}}]\label{Thm:locally-finite-KG}
Let $\class{K}$ be a subvariety of $\class{KG}$. The following conditions are equivalent:
\benroman
\item $\class{K}$ is locally finite.
\item $\class{K}$ excludes an algebra of the form $\A + \bs{2}$ where $\A$ is a finite FSI member of $\HHH(\bs{RN})$. 
\item $\class{K}$ excludes the algebra $\bs{RN} + \bs{2}$.
\eroman
\end{Theorem}


We shall rely also on some observations about the topological duals of the members of $\class{KG}$ (Corollary \ref{Cor:dual-of-KG}). Let $n \in \omega$. A poset $\boldsymbol{X} = \langle X; \leq \rangle$ is said to have \emph{incomparability degree} at most $n$ if for every $x \in X$, the set ${\uparrow} x$ does not contain any point which is incomparable with $n+1$ elements. Clearly, posets of incomparability degree at most $n$ also have width at most $n + 1$, but the converse need not be true in general (since elements incomparable with a given element may be comparable with each other). 
A Heyting algebra $\A$ has \emph{incomparability degree} at most $n$, when this is the case for the poset underlying its dual space $\A_{\ast}$. 
We denote by $\class{ID}_{n}$ the class of all Heyting algebras of incomparability degree at most $n$. We shall see that $\class{ID}_{n}$ is 
a variety.

Let $n \in \omega$. Consider a set of variables $Z_{n} = \{ y_{1}, \dots, y_{n+1} \}$. We let $\mathbb{Z}_{n,1}, \dots, \mathbb{Z}_{n, k_n}$ be a fixed enumeration of all possible posets with universe $Z_{n}$. For each such $\mathbb{Z}_{n, k} = \langle Z_n,\leq_{k} \rangle$, with $k \leq k_n$, define the formulas
\[
\psi_{n, k} \coloneqq \bigvee_{i=1}^{n+1}(y_{i} \to (x \lor \bigvee_{\substack{j \colon y_{i} \nleq_{k} y_{j}}} y_{j})).
\]
In the above display we assume that the disjunction of an empty family is the symbol $0$. Moreover, we set
\[
\delta_{n, k} \coloneqq 
\psi_{n, k} \lor  (x \to \bigvee_{i=1}^{n+1} y_{i})
\]
and
\[
\Sigma_{n} \coloneqq \{ \delta_{n, k} \thickapprox 1 \colon k = 1, \dots, k_n \}.
\]

\begin{Theorem}
For every $n \in \omega$, the class $\class{ID}_{n}$ of Heyting algebras with incomparability degree at most $n$ is axiomatized by the set of equations $\Sigma_{n}$. As a consequence, $\class{ID}_{n}$ is a variety.
\end{Theorem}

\begin{proof}
First we show that for every Heyting algebra $\A \notin \class{ID}_n$, we have $\A \not\models \Sigma_{n}$. 
Note that we need only exhibit the failure of some equation of $\Sigma_{n}$ in some homomorphic image of $\A$. Since $\A \notin \class{ID}_n$, there is an $x \in \A_{\ast}$ such that ${\uparrow} x$ contains a point which is incomparable with $n+1$ points. Thus, owing to the previous remark, and (\ref{Lem:correspondences:FSI}) and (\ref{Lem:correspondences:FactorAndSubspaces}) of Lemma~\ref{Lem:correspondences}, we may, without loss of generality, suppose that $\A$ is FSI, otherwise we replace 
$\A$ with its FSI homomorphic image 
whose dual is isomorphic
to the subspace ${\uparrow} x$.

Since $\A \notin \class{ID}_n$, there are distinct $F, G_{1}, \dots, G_{n+1} \in \Pr \A$ such that $F$ is incomparable with 
each of 
$G_{1}, \dots, G_{n+1}$. Then for every $i \leq n+1$ we can choose an element $a_{i} \in F \smallsetminus G_{i}$. We set
\[
\hat{a} \coloneqq a_{1} \land \dots \land a_{n+1}.
\]
Observe that
\begin{equation}\label{Eq:the-element-a}
\hat{a} \in F \smallsetminus (G_{1} \cup \dots \cup G_{n+1}).
\end{equation}

Given $i \leq n+1$, we can choose $b_{i} \in G_{i} \smallsetminus F$. Moreover, for every $j \leq n+1$ such that $G_{i} \nsubseteq G_{j}$, we choose $b_{i_j} \in G_{i} \smallsetminus G_{j}$ and set
\[
b_{i}' \coloneqq b_{i} \land \bigwedge_{\substack{j \colon G_{i} \nsubseteq G_{j}}} b_{i_{j}}.
\]
Finally for every $j \leq n+1$, we define
\[
\hat{b}_{j} \coloneqq \bigwedge_{\substack{i \colon G_{i} \subseteq G_{j}}} b_{i}'.
\]
Observe that for every $i \leq n+1$,
\begin{equation}\label{Eq:the-elements-b}
\hat{b}_{i} \in G_{i} \smallsetminus (F \cup \bigcup_{\substack{j \colon G_{i}\nsubseteq G_{j}}} G_{j}).
\end{equation}

From (\ref{Eq:the-element-a}, \ref{Eq:the-elements-b}) and the fact that $G_{1}, \dots, G_{n+1}$ are different, we deduce 
that the elements $\hat{a}, \hat{b}_{1}, \dots, \hat{b}_{n+1}$ are different one from the other. Then consider the subposet $\mathbb{P}$ of $\A$ with universe $\{ \hat{b}_{1}, \dots, \hat{b}_{n+1} \}$. Clearly there is a $k \leq k_n$ such that $\mathbb{Z}_{n, k}$ is isomorphic to $\mathbb{P}$ under the map
$
y_{i} \mapsto \hat{b}_{i} \;(i \leq n+1). 
$
We leave it for the reader to verify that $\delta_{n,k}^{\A}(\hat{a}, \hat{b}_1, \dots, \hat{b}_{n+1}) \neq 1$. It may be helpful to notice that for every $i, j  \leq n+1$,
\begin{equation}\label{Eq:the-elements-c}
\hat{b}_{i} \in G_{j} \Longleftrightarrow \hat{b}_{j} \leq \hat{b}_{i}.
\end{equation}

Conversely, we show that for every Heyting algebra $\A$ such that $\A \nvDash \Sigma_{n}$, we get $\A \notin \class{InD}_n$.
There exists $k \leq k_n$ such that the equation $\delta_{n, k} \thickapprox 1$ in $\Sigma_{n}$ is not valid in $\A$, so 
there are $a, b_{1}, \dots, b_{n+1}$ such that
%
\[
\bigvee_{i \leq n+1}({b}_{i} \to ({a} \lor \bigvee_{\substack{j \colon y_{i} \nleq_{k} y_{j}}} {b}_{j})) \lor ({a} \to \bigvee_{i \leq n+1} {b}_{i}) \ne 1.
\]
Therefore,
for every $i \leq n+1$, 
\[
{b}_{i} \nleq {a} \lor \bigvee_{\substack{j \colon y_{i} \nleq_{k} y_{j}}} {b}_{j}, \quad\text{ and }\quad
{a}  \nleq \bigvee_{j \leq n+1} {b}_{j}.
\]
By the \emph{prime filter lemma} \cite[Thm.~4.1]{Ber12} for distributive lattices, there are
prime filters $F, G_{1}, \dots, G_{n+1}$ of $\A$ 
such that
\[
{a} \in F \text{ and }{b}_{1}, \dots, {b}_{n+1} \notin F
\]
and for every $i \leq n+1$,
\[
{b}_{i} \in G_{i} \text{ and }{a} \cup \{ {b}_{j} \colon y_{i} \nleq_{k} y_{j} \} \subseteq A \smallsetminus G_{i}. 
\]
It follows easily from these properties and the pairwise distinctness of $y_1, \dots, y_{n+1}$
that $G_{1}, \dots, G_{n+1}$ are pairwise different, and that $F$ is incomparable with every $G_{i}$ for every $i \in \{1, \dots, n+1 \}$ 
Therefore 
$\A \notin \class{ID}_{n}$ 
as required.
\end{proof}

\begin{Corollary}\label{Cor:dual-of-KG}
$\class{KG} \subseteq \class{ID}_{2} \cap \class{W}_{2}$.
\end{Corollary}

\begin{proof}
The generators of $\class{KG}$ in (\ref{Eq:generators-KG}) are easily seen to belong to $\class{ID}_{2} \cap \class{W}_{2}$. Since both $\class{ID}_{2}$ and $\class{W}_{2}$ are varieties, this implies that $\class{KG} \subseteq \class{ID}_{2} \cap \class{W}_{2}$.
\end{proof}

\section{A continuum of failures of the ES property}

In this section we show that there is a continuum of locally finite subvarieties of $\VVV(\bs{RN})$ lacking the ES property (Theorem \ref{Thm:continuum}). By Corollary \ref{Cor:dual-of-KG}, there is also such a continuum among varieties with width at most $2$. 

Recall, from the end of Section~\ref{sec:sums}, that $\bs{D}_{2}$ is the two-element discrete poset, equipped with the discrete topology. 

\begin{Lemma}\label{Lem:D2-in-RN}
$(\bs{D}_{2}^{\infty})^{\ast} \in \VVV(\bs{RN})$.
\end{Lemma}

\begin{proof}
Recall that every algebra embeds into a ultraproduct of its finitely generated subalgebras \cite[Thm.\ 2.14, Ch.\ V]{BS81}. Now, observe that the finitely generated subalgebras of $(\bs{D}_{2}^{\infty})^{\ast}$ coincide with finite sums, where each summand is either the two-element or the four-element Boolean algebra. 
It is therefore not hard to see, when considering the subalgebra of encircled elements in the figure of $\bs{RN}$ below
on page~\pageref{fig:RN-subalgebra} (Figure~\ref{fig:RN-subalgebra}), 
that every finitely generated subalgebra of $(\bs{D}_{2}^{\infty})^{\ast}$ belongs to $\HHH\SSS(\bs{RN})$.
As a consequence, we obtain that $(\bs{D}_{2}^{\infty})^{\ast} \in \SSS\PPU\HHH\SSS(\bs{RN}) \subseteq \VVV(\bs{RN})$.
\end{proof}

\begin{Lemma}\label{Lem:excluding-D2}
Let $\class{K}$ be a subvariety of $\class{ID}_{2} \cap \class{W}_{2}$. If $\bs{D}_{2}^{\infty} \in \class{K}_{\ast}$, then $\class{K}$ lacks the ES property.
\end{Lemma}

\begin{proof}

Suppose with a view to obtaining a contradiction 
that there is a variety $\class{K} \subseteq \class{ID}_{2} \cap \class{W}_{2}$ with the ES property and such that $\bs{D}_{2}^{\infty} \in \class{K}_{\ast}$. 
Let $R$ be the equivalence relation on $\bs{D}_{2}^{\infty}$ whose corresponding partition is depicted in the diagram below:
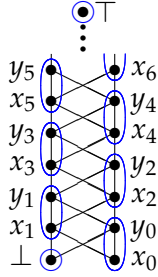
\begin{center}
\begin{picture}(48,112)

\newsavebox{\Botie}
\savebox{\Botie}(24,12){%
\put(0,0){\line(2,1){24}}
\put(24,12){\circle*{4}}
\put(0,0){\line(0,1){12}}
\put(0,12){\circle*{4}}
\put(0,12){\line(2,-1){24}}
\put(24,0){\line(0,1){12}}
}
\put(12,12){\circle*{4}}
\put(36,12){\circle*{4}}
\put(12,12){\usebox{\Botie}}
\put(12,24){\usebox{\Botie}}
\put(12,36){\usebox{\Botie}}
\put(12,48){\usebox{\Botie}}
\put(12,60){\usebox{\Botie}}
\put(12,72){\usebox{\Botie}}
\put(12,84){\line(0,1){6}}
\put(36,84){\line(0,1){6}}

{\color{blue}

\newsavebox{\EqR}
\savebox{\EqR}(8,20){
\qbezier(4,0)(0,0)(0,10)
\qbezier(0,10)(0,20)(4,20)
\qbezier(4,0)(8,0)(8,10)
\qbezier(8,10)(8,20)(4,20)%
}

\put(32,8){\usebox{\EqR}}
\put(32,32){\usebox{\EqR}}
\put(32,56){\usebox{\EqR}}
\put(8,20){\usebox{\EqR}}
\put(8,44){\usebox{\EqR}}
\put(8,68){\usebox{\EqR}}

\qbezier(36,80)(32,80)(32,90)
\qbezier(36,80)(40,80)(40,90)

\put(12,12){\circle{8}}
\put(24,106){\circle{8}}

}

\put(24,92){\circle*{2}}
\put(24,96){\circle*{2}}
\put(24,100){\circle*{2}}

\put(24,106){\circle*{4}}

\put(-4,10){$\bot$}
\put(-4,22){$x_{1}$}
\put(-4,34){$y_{1}$}
\put(-4,46){$x_{3}$}
\put(-4,58){$y_{3}$}
\put(-4,70){$x_{5}$}
\put(-4,82){$y_{5}$}

\put(42,10){$x_{0}$}
\put(42,22){$y_{0}$}
\put(42,34){$x_{2}$}
\put(42,46){$y_{2}$}
\put(42,58){$x_{4}$}
\put(42,70){$y_{4}$}
\put(42,82){$x_{6}$}

\put(28,102){$\top$}
\end{picture}
\end{center}
An argument analogous to the one detailed for Lemma \ref{Lem:correct-partition} shows that $R$ is a correct partition on $\bs{D}_{2}^{\infty}$. Since $\class{K}$ has the ES property, we can apply Lemma \ref{Lem:epic-subalgebras-dual}, so 
there exist $\boldsymbol{Y} \in \class{K}_{\ast}$ and a pair of different Esakia morphisms $f, g \colon \boldsymbol{Y} \to \bs{D}_{2}^{\infty}$ such that  $\langle f(y), g(y) \rangle \in R$ for every $y \in Y$.

Since $f \ne g$, there exists $\bot \in Y$ such that $f(\bot) \ne g(\bot)$. Together with the fact that $\langle f(\bot), g(\bot) \rangle \in R$, this implies that $\{ f(\bot), g(\bot) \} = \{ x_{n}, y_{n} \}$ for some $n \in \omega$. We can assume without loss of generality that $f(\bot) = x_{n}$ and $g(\bot) = y_{n}$. Moreover, as in the proof of Theorem \ref{Thm:Es-bounded-width}, we can assume that $Y = {\uparrow} \bot$. Observe that 
$\bs{D}_{2}^{\infty}$ and $\bs{Y}$
have width at most $2$.
Moreover, $\bs{Y}$ and $\bs{D}_{2}^{\infty}$, and the Esakia morphism $f \colon \boldsymbol{Y} \to \bs{D}_{2}^{\infty}$ satisfy the assumptions of Lemma \ref{Lem:trick-width}. Therefore, $\boldsymbol{Y}$ has a subposet $\langle Z; \leq^{\boldsymbol{Y}}\rangle$ such that the restriction
\[
f \colon \langle Z; \leq^{\boldsymbol{Y}}\rangle \to \langle {\uparrow} f(\bot)^{\top}; \leq^{\bs{D}_{2}^{\infty}}\rangle
\]
is a poset isomorphism. For the sake of simplicity, we denote the elements of $Z$ exactly as their alter egos in ${\uparrow} f(\bot)^{\top}$. Under this convention, 
\[
Z = \{ x_{n+ m} \colon m \in \omega \} \cup \{ y_{n + m} \colon m \in \omega \}
\]
and 
\[
f(x_{i}) = x_{i} \text{ and }f(y_{i}) = g(x_{i}) = g(y_{i}) = y_{i}
\]
for every $x_{i}, y_{i} \in Z$ (see the proof of Theorem \ref{Thm:Es-bounded-width}, if necessary).

We shall 
now investigate 
the structure of $\bs{Y}$ to produce a sequence of elements $z_{n+2},z_{n+3},\,\ldots \in Y \smallsetminus Z$ and describe how they are ordered with respect to the elements of $Z$. 
First, observe that $g(y_{n}) = y_{n} \leq^{\bs{D}_{2}^{\infty}} x_{n+2}$. Since $g$ is an Esakia morphism, there is an element $z_{n+2} \in Y$ such that $y_{n} \leq^{\bs{Y}} z_{n+2}$ and $g(z_{n+2}) = x_{n+2}$.
Let us describe the structure of the poset $\langle Z \cup \{ z_{n+2} \}; \leq^{\bs{Y}} \rangle$. First observe that $z_{n+2}$ is incomparable with $x_{n+1}$ and $y_{n+1}$ with respect to $\leq^{\bs{Y}}$, since $g$ is order-preserving and $g(z_{n+2}) = x_{n+2}$ is incomparable with $g(x_{n+1}) = g(y_{n+1}) = y_{n+1}$ in $\bs{D}_{2}^{\infty}$. Moreover, $z_{n+2} \leq^{\bs{Y}} x_{n+2}$. To prove this, observe that $x_{n+2}$ and $y_{n+1}$ are incomparable in $\bs{Y}$. Since $\bs{Y}$ has width at most $2$, this implies that $z_{n+2}$ must be comparable with one of them. Since $z_{n+2}$ is incomparable with $y_{n+1}$, if follows that $z_{n+2}$ is comparable with $x_{n+2}$. Keeping in mind that $g(x_{n+2}) = y_{n+2} \nleq^{\bs{D}_{2}^{\infty}} x_{n+2} = g(z_{n+2})$ and that $g$ is order-preserving, we obtain 
$x_{n+2} \nleq^{\bs{Y}} z_{n+2}$. As a consequence, we conclude that $z_{n+2} <^{\bs{Y}} x_{n+2}$ as desired. Summing up, the structure of $\langle Z \cup \{ z_{n+2} \}; \leq^{\bs{Y}} \rangle$ is described exactly by the following picture:

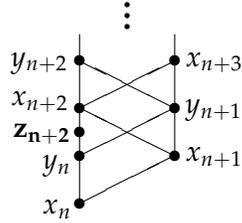
\begin{center}
\begin{picture}(100,85)
\put(26,4){\circle*{4}}
\put(26,4){\line(0,1){64}}
\put(26,4){\line(2,1){36}}
\put(62,22){\circle*{4}}
\put(62,22){\line(0,1){46}}
\put(62,22){\line(-2,1){36}}
\put(26,22){\circle*{4}}
\put(26,22){\line(2,1){36}}
\put(62,40){\circle*{4}}
\put(62,40){\line(-2,1){36}}
\put(26,40){\circle*{4}}
\put(26,40){\line(2,1){36}}
\put(62,58){\circle*{4}}
\put(26,58){\circle*{4}}

\put(26,31){\circle*{4}}

\put(44,71){\circle*{2}}
\put(44,75){\circle*{2}}
\put(44,79){\circle*{2}}


\put(12,0){$x_{n}$}
\put(12,17){$y_{n}$}
\put(1,29){$\mathbf{z_{n+2}}$}
\put(0,41){$x_{n+2}$}
\put(0,56){$y_{n+2}$}
\put(66,20){$x_{n+1}$}
\put(66,38){$y_{n+1}$}
\put(66,56){$x_{n+3}$}

\end{picture}
\end{center}

Now, observe that $g(z_{n+2}) = x_{n+2} \leq^{\bs{D}_{2}^{\infty}} x_{n+3}$. Since $g$ is an Esakia morphism, there is an element $z_{n+3} \in Y$ with $z_{n+2} \leq^{\bs{Y}} z_{n+3}$ such that $g(z_{n+3}) = x_{n+3}$. We can replicate the previous argument, used to describe the structure of the poset $\langle Z \cup \{ z_{n+2} \}; \leq^{\bs{Y}} \rangle$, to show that $z_{n+3}$ is incomparable with $x_{n+2}$ and $y_{n+2}$, and that $z_{n+3} <^{\bs{Y}} x_{x+3}$. Then, as in the previous argument, $y_{n+1} <^{\bs{Y}} z_{n+3}$.
Iterating this process we construct a series of elements $\{ z_{n + m} \colon 2 \leq m \in \omega \} \subseteq Y$ such that $g(z_{i}) = x_{i}$, for all $i \geq 2$. 
The structure of the poset $\mathbb{Z}' \coloneqq \langle Z \cup \{ z_{n+m} \colon 2 \leq m \in \omega \}; \leq^{\bs{Y}} \rangle$ is as 
depicted below:
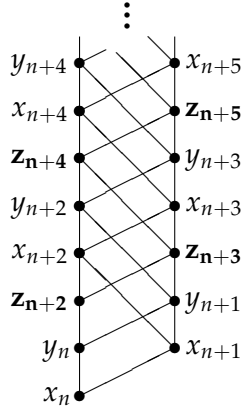
\begin{center}
\begin{picture}(100,160)
\put(26,4){\line(0,1){136}}
\put(62,22){\line(0,1){118}}
\put(26,4){\circle*{4}}
\put(26,4){\line(2,1){36}}
\put(62,22){\circle*{4}}
\put(62,22){\line(-1,1){36}}
\put(26,22){\circle*{4}}
\put(26,22){\line(2,1){36}}
\put(62,40){\circle*{4}}
\put(62,40){\line(-1,1){36}}
\put(26,40){\circle*{4}}
\put(26,40){\line(2,1){36}}
\put(62,58){\circle*{4}}
\put(26,58){\circle*{4}}

\put(62,58){\line(-1,1){36}}
\put(26,58){\circle*{4}}
\put(26,58){\line(2,1){36}}
\put(62,76){\circle*{4}}
\put(26,76){\circle*{4}}

\put(62,76){\line(-1,1){36}}
\put(26,76){\circle*{4}}
\put(26,76){\line(2,1){36}}
\put(62,94){\circle*{4}}
\put(26,94){\circle*{4}}

\put(62,94){\line(-1,1){36}}
\put(26,94){\circle*{4}}
\put(26,94){\line(2,1){36}}
\put(62,112){\circle*{4}}
\put(26,112){\circle*{4}}

\put(62,112){\line(-1,1){22}}
\put(26,112){\circle*{4}}
\put(26,112){\line(2,1){36}}
\put(62,130){\circle*{4}}
\put(26,130){\circle*{4}}

\put(62,130){\line(-1,1){10}}
\put(26,130){\line(2,1){10}}

\put(44,144){\circle*{2}}
\put(44,148){\circle*{2}}
\put(44,152){\circle*{2}}


\put(12,2){$x_{n}$}
\put(12,20){$y_{n}$}
\put(0,38){$\mathbf{z_{n+2}}$}
\put(0,56){$x_{n+2}$}
\put(0,74){$y_{n+2}$}
\put(0,92){$\mathbf{z_{n+4}}$}
\put(0,110){$x_{n+4}$}
\put(0,128){$y_{n+4}$}
\put(66,20){$x_{n+1}$}
\put(66,38){$y_{n+1}$}
\put(66,56){$\mathbf{z_{n+3}}$}
\put(66,74){$x_{n+3}$}
\put(66,92){$y_{n+3}$}
\put(66,110){$\mathbf{z_{n+5}}$}
\put(66,128){$x_{n+5}$}

\end{picture}
\end{center}

We claim that for every $a \in Y$ such that $x_{n+2} \leq^{\bs{Y}} a$, either $a \in Z'$ or $b \leq^{\bs{Y}} a$ for every $b \in Z'$. To prove this, consider $a \in Y$ such that $x_{n+2} \leq^{\bs{Y}} a$ and $a \notin Z'$. It will be enough to show that $x_{n+m} \leq^{\bs{Y}}a$ for $2 < m \in \omega$. Suppose, with a view to obtaining a contradiction, 
that there is a smallest integer $m > 2$ such that $x_{n+m} \nleq^{\bs{Y}} a$. Looking at the figure above, it is easy to see that every point in $Z' \smallsetminus \{ x_{n}, y_{n} \}$ is incomparable with two elements in ${\uparrow} x_{n}$. Since $\bs{Y}$ has incomparability degree at most $2$ by Corollary \ref{Cor:dual-of-KG}, it follows that every element in ${\uparrow} x_{n} \smallsetminus Z'$ is comparable with all the elements of $Z' \smallsetminus \{ x_{n}, y_{n} \}$. We shall make extensive use of this observation. First recall that $x_{n+m} \nleq^{\bs{Y}} a$. As $a$ is comparable with $x_{n+m}$, this implies that $a <^{\bs{Y}} x_{n+m}$. Moreover, $a$ is comparable with $y_{n+m-1}$. Since $y_{n+m-1} \nleq^{\bs{Y}} x_{n+m}$ and $a<^{\bs{Y}} x_{n+m}$, it follows that $a <^{\bs{Y}} y_{n+m-1}$. Now, $a$ is comparable with $z_{n+m}$. If $z_{n+m} \leq^{\bs{Y}} a$, then $z_{n+m} \leq^{\bs{Y}} y_{n+m-1}$, which is false. 
Thus, $a <^{\bs{Y}} z_{n+m}$. By minimality of $m$ we have $x_{n+m-1} \leq^{\bs{Y}}a$. This yields that $x_{n+m-1}\leq^{\bs{Y}} z_{n+m}$. This contradiction establishes the claim.

From the definition of an Esakia space we know that the upset ${\uparrow} x_{n+2}$ in $\bs{Y}$ is closed and, therefore, an E-subspace of $\boldsymbol{Y}$. Moreover, ${\uparrow} x_{n+2} \in \class{K}_{\ast}$, since $\class{K}$ is a variety. Consider the equivalence relation $S$ on ${\uparrow} x_{n+2}$ defined as follows: for every $a, b \in {\uparrow} x_{n+2}$,
\begin{align*}
\langle a, b \rangle \in S \Longleftrightarrow& \text{ either }a = b \text{ or }a, b  \notin Z'.
\end{align*}
We shall prove that $S$ is a correct partition on ${\uparrow} x_{n+2}$. To this end, observe that from the claim it follows that $S$ satisfies condition (i) in the definition of a correct partition. In order to prove condition (ii), consider $a, b \in {\uparrow} x_{n+2}$ such that $\langle a, b \rangle \notin S$.
We can assume without loss of generality that $b \in Z'$. If $b \in \{x_{n+2}, x_{n+2}\}$, 
let $b' = z_{n+4}$; 
otherwise, 
let $b' = b$. Let $c$ be the minimum element of ${\uparrow} x_{n+2}$ that is incomparable with $b'$. By the Priestley separation axiom, since $c \nleq^{\bs{Y}}b'$, there is a clopen upset $U$ such that $c \in U$ and $b' \notin U$. Looking at the above picture, it is easy to see that the only upset missing $b'$ is ${\uparrow} c$, therefore $U = {\uparrow} c$ and $U^{c} = {\downarrow} b'$. In particular, 
$b \in {\downarrow} b' = U^{c}$. By the claim above, $a \in {\uparrow} c = U$. The fact that $U$ and $U^{c}$ are unions of equivalence classes of $S$ follows from the definition of $S$. This establishes condition (ii) and, therefore, that $S$ is a correct partition on ${\uparrow} x_{n+2}$.

Then let $\bs{W}$ be the Esakia space $({\uparrow} x_{n+2} )/ S$. Observe that $\bs{W} \in \class{K}_{\ast}$, since $\class{K}$ is 
closed under homomorphic images.
Moreover, observe that the poset underlying $\bs{W}$ is isomorphic to $Z' \cap {\uparrow} x_{n+2}$ plus a fresh top element.\footnote{The existence of the top element of $\bs{W}$ follows from the fact that, in Esakia spaces, suprema of chains exist, see Lemma \ref{Lem:correspondences}(\ref{Lem:correspondences:Chains}).}
Now, consider the equivalence relation $T$ on $W$ 
whose corresponding partition is 
depicted below:
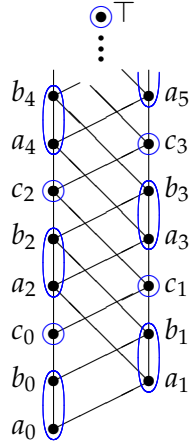
\begin{center}
\begin{picture}(100,165)
\put(26,4){\line(0,1){136}}
\put(62,22){\line(0,1){118}}
\put(26,4){\circle*{4}}
\put(26,4){\line(2,1){36}}
\put(62,22){\circle*{4}}
\put(62,22){\line(-1,1){36}}
\put(26,22){\circle*{4}}
\put(26,22){\line(2,1){36}}
\put(62,40){\circle*{4}}
\put(62,40){\line(-1,1){36}}
\put(26,40){\circle*{4}}
\put(26,40){\line(2,1){36}}
\put(62,58){\circle*{4}}
\put(26,58){\circle*{4}}

\put(62,58){\line(-1,1){36}}
\put(26,58){\circle*{4}}
\put(26,58){\line(2,1){36}}
\put(62,76){\circle*{4}}
\put(26,76){\circle*{4}}

\put(62,76){\line(-1,1){36}}
\put(26,76){\circle*{4}}
\put(26,76){\line(2,1){36}}
\put(62,94){\circle*{4}}
\put(26,94){\circle*{4}}

\put(62,94){\line(-1,1){36}}
\put(26,94){\circle*{4}}
\put(26,94){\line(2,1){36}}
\put(62,112){\circle*{4}}
\put(26,112){\circle*{4}}

\put(62,112){\line(-1,1){22}}
\put(26,112){\circle*{4}}
\put(26,112){\line(2,1){36}}
\put(62,130){\circle*{4}}
\put(26,130){\circle*{4}}

\put(62,130){\line(-1,1){10}}
\put(26,130){\line(2,1){10}}

{\color{blue}

\newsavebox{\EqRa}
\savebox{\EqRa}(8,26){
\qbezier(4,0)(0,0)(0,13)
\qbezier(0,13)(0,26)(4,26)
\qbezier(4,0)(8,0)(8,13)
\qbezier(8,13)(8,26)(4,26)%
}
\put(22,0){\usebox{\EqRa}}
\put(22,54){\usebox{\EqRa}}
\put(22,108){\usebox{\EqRa}}
\put(58,18){\usebox{\EqRa}}
\put(58,72){\usebox{\EqRa}}

\qbezier(62,126)(58,126)(58,139)
\qbezier(62,126)(66,126)(66,139)

\put(26,40){\circle{8}}
\put(26,94){\circle{8}}
\put(62,58){\circle{8}}
\put(62,112){\circle{8}}
\put(44,160){\circle{8}}

}

\put(44,144){\circle*{2}}
\put(44,148){\circle*{2}}
\put(44,152){\circle*{2}}

\put(44,160){\circle*{4}}

\put(10,2){$a_{0}$}
\put(10,20){$b_{0}$}
\put(10,38){$c_{0}$}
\put(10,56){$a_{2}$}
\put(10,74){$b_{2}$}
\put(10,92){$c_{2}$}
\put(10,110){$a_{4}$}
\put(10,128){$b_{4}$}
\put(68,20){$a_{1}$}
\put(68,38){$b_{1}$}
\put(68,56){$c_{1}$}
\put(68,74){$a_{3}$}
\put(68,92){$b_{3}$}
\put(68,110){$c_{3}$}
\put(68,128){$a_{5}$}

\put(48,158){$\top$}
\end{picture}
\end{center}
An argument, similar to the one detailed in the case of $S$, shows that the relation $T$ is a correct partition on $\bs{W}$, except that in this case we let $b$ be such that $a \nleq b$, and let $b'$ be $c_0$ 
if $b \in \{a_0,b_0\}$ and the maximum of the equivalence class $b/T$ otherwise.

Since $\class{K}$ has the ES property, we can apply Lemma \ref{Lem:epic-subalgebras-dual}, so 
there exist $\boldsymbol{V} \in \class{K}_{\ast}$ and a pair of different Esakia morphisms $f, g \colon \boldsymbol{V} \to \bs{W}$ such that  $\langle f(v), g(v) \rangle \in T$ for every $v \in V$. As above, since $f \ne g$, there are $\bot \in V$ and $n \in \omega$ such that $\{ f(\bot), g(\bot) \} = \{ a_{n}, b_{n} \}$. We can assume without loss of generality that $f(\bot) = a_{n}$ and $g(\bot) = b_{n}$, and that $V = {\uparrow} \bot$. 
Moreover, we can find a subposet $\langle Q; \leq^{\bs{V}}\rangle$ of $\bs{V}$ such that the restriction
\[
f \colon \langle Q; \leq^{\boldsymbol{V}}\rangle \to \langle {\uparrow} f(\bot)^{\top}; \leq^{\bs{W}}\rangle
\]
is a poset isomorphism. We denote the elements of $Q$ exactly as their alter egos in ${\uparrow} f(\bot)^{\top}$. Under this convention, 
\[
Q = \{ a_{n+ m} \colon m \in \omega \} \cup \{ b_{n + m} \colon m \in \omega \}\cup \{ c_{n + m} \colon m \in \omega \}
\]
and 
for every $a_{i}, b_{i}, c_{i} \in Q$,
\[
f(a_{i}) = a_{i} \text{ and }f(b_{i}) = g(a_{i}) = g(b_{i}) = b_{i} \text{ and }f(c_{i}) = g(c_{i}) = c_{i}.
\]
Observe that $g(b_{n}) = b_{n} \leq^{\bs{W}} a_{n+2}$. Since $g$ is an Esakia morphism, there exists 
$v \in V$ such that $b_{n} \leq^{\bs{V}} v$ and $g(v) = a_{n+2}$. So,
we have that $\{ g(v), g(a_{n+2}), g(b_{n+2}) \} = \{ a_{n+2}, b_{n+2} \}$ and $g(c_{n+1}) = c_{n+1}$, therefore, the elements $g(v), g(a_{n+2}), g(b_{n+2})$ are incomparable with $g(c_{n+1})$ in $\bs{W}$. Since $g$ is order-preserving, 
$v, a_{n+2}, b_{n+2}$ are incomparable with $c_{n+1}$ in $\bs{V}$. 
Because $a_{n+2} \ne b_{n+2}$ and 
$\bs{V}$ has incomparability degree $\leq 2$ by Corollary \ref{Cor:dual-of-KG}, we conclude that either $v = a_{n+2}$ or $v= b_{n+2}$. Observe that if $v= a_{n+2}$, then $a_{n+2} = g(v) = g(a_{n+2}) = b_{n+2}$, which is false. A similar argument rules out the case where $v= b_{n+2}$. Hence we have reached a contradiction, as desired.
\end{proof}

\begin{Corollary}\label{Cor:RN-not-ES}
The members of the interval $[ \VVV((\bs{D}_{2}^{\infty})^{\ast}), \class{ID}_{2} \cap \class{W}_{2}]$ of the subvariety lattice of $\class{HA}$ lack the ES property. In particular, $\VVV(\bs{RN})$ lacks the ES property.
\end{Corollary}
\begin{proof}
From Lemma \ref{Lem:D2-in-RN} and Corollary \ref{Cor:dual-of-KG} we know that $\VVV(\bs{RN})$ belongs to the interval $[ \VVV((\bs{D}_{2}^{\infty})^{\ast}), \class{ID}_{2} \cap \class{W}_{2}]$. The fact that $\VVV(\bs{RN})$ lacks the ES property is then an immediate consequence of Lemma \ref{Lem:excluding-D2}.
\end{proof}

As already mentioned, up to now the only known example of a variety of Heyting algebras without the ES property was precisely $\VVV((\bs{D}_{2}^{\infty})^{\ast})$. This example is now subsumed by the above corollary.

\begin{Theorem}\label{Thm:continuum}
There is a continuum of locally finite subvarieties of $\VVV(\bs{RN})$ without the ES property. 
\end{Theorem}

\begin{proof}
Let $\C$ be the three-element Heyting algebra. Moreover, define $\A \coloneqq \bs{2} \times \C$. For $2 \leq n \in \omega$, let $\B_{n}$ be the algebra $\bs{2} + \A + \C_{1} + \dots +  \C_{n-2}$, where $\C_{1}, \dots, \C_{n-2}$ are copies of the four-element Boolean algebra ${\bs{D}_{2}}^{\ast}$. 
The algebra $\B_n$ (depicted below) 
belongs to $\HHH\SSS(\bs{RN})$. 
\begin{center}
\begin{picture}(100,165)
\put(36,4){\circle*{4}}
\put(36,4){\line(1,1){12}}
\put(36,4){\line(-1,1){12}}
\put(48,16){\circle*{4}}
\put(24,16){\circle*{4}}
\put(36,28){\line(1,-1){12}}
\put(36,28){\line(-1,-1){12}}
\put(36,28){\circle*{4}}

\put(36,36){\circle*{2}}
\put(36,40){\circle*{2}}
\put(36,44){\circle*{2}}

\put(36,52){\circle*{4}}
\put(36,52){\line(1,1){12}}
\put(36,52){\line(-1,1){12}}
\put(48,64){\circle*{4}}
\put(24,64){\circle*{4}}
\put(36,76){\line(1,-1){12}}
\put(36,76){\line(-1,-1){12}}
\put(36,76){\circle*{4}}

\put(36,76){\line(1,1){12}}
\put(36,76){\line(-1,1){12}}
\put(48,88){\circle*{4}}
\put(24,88){\circle*{4}}
\put(36,100){\line(1,-1){12}}
\put(36,100){\line(-1,-1){12}}
\put(36,100){\circle*{4}}

\put(36,100){\line(1,1){12}}
\put(36,100){\line(-1,1){12}}
\put(48,112){\circle*{4}}
\put(24,112){\circle*{4}}
\put(36,124){\line(1,-1){12}}
\put(36,124){\line(-1,-1){12}}
\put(36,124){\circle*{4}}

\put(24,112){\line(-1,1){12}}
\put(12,124){\circle*{4}}
\put(24,136){\line(1,-1){12}}
\put(24,136){\line(-1,-1){12}}
\put(24,136){\circle*{4}}

\put(24,148){\line(0,-1){12}}
\put(24,148){\circle*{4}}

\put(54,30){\scalebox{1}[8]{$\}$}}
\put(64,54){$n-2$}
\put(64,43){copies}
\put(64,30){of ${\bs{D}_{2}}^{\ast}$}
\end{picture}
\end{center}
Define $F \coloneqq \{ \B_{n} \colon 2 \leq n \in \omega \}$. In \cite[Lem.\ 5.38(5), Thm.\ 5.39(1)]{BBdeJ08} it is shown that $\VVV(S) \ne \VVV(T)$, for every pair of different 
subsets $S, T \subseteq F$.\footnote{The proof relies on a general method introduced by Jankov for the case of varieties of Heyting algebras \cite{Jankov68}, and subsequently extended to all varieties with equationally definable principal congruences by Blok and Pigozzi \cite{BP82} (see also \cite{Jon95}).}

We claim that $\VVV(S, (\bs{D}_{2}^{\infty})^{\ast}) \ne \VVV(T, (\bs{D}_{2}^{\infty})^{\ast})$, for every pair of different subsets $S, T \subseteq F$. To prove this, consider two different $S, T \subseteq F$. Since $\VVV(S) \ne \VVV(T)$, we can assume without loss of generality that $\B_{n} \in \VVV(S) \smallsetminus \VVV(T)$ whenever $n \geq 2$. Suppose with a view to obtaining a contradiction 
that $\VVV(S, (\bs{D}_{2}^{\infty})^{\ast}) = \VVV(T, (\bs{D}_{2}^{\infty})^{\ast})$. In particular, $\B_{n} \in \VVV(T, (\bs{D}_{2}^{\infty})^{\ast})_{\textup{FSI}}$. Now, from J\'onsson's lemma it follows that $\VVV(T, (\bs{D}_{2}^{\infty})^{\ast})_{\textup{FSI}} = \VVV(T)_{\textup{FSI}} \cup \VVV((\bs{D}_{2}^{\infty})^{\ast})_{\textup{FSI}}$. Since $\B_{n} \notin \VVV(T)_{\textup{FSI}}$, we have that $\B_{n} \in \VVV((\bs{D}_{2}^{\infty})^{\ast})$. Now, observe that $(\bs{D}_{2}^{\infty})^{\ast} \in \class{ID}_{1}$. As a consequence, $\B_{n} \in \class{ID}_{1}$. But this is easily seen to be false. 
So, we have reached a contradiction, thus establishing the claim.

From the claim it follows that the set
\[
G \coloneqq \{ \VVV(T, \D_{2}^{\infty}) \colon T \subseteq F \}
\]
has the cardinality of the continuum. Consider $T \subseteq F$. It remains only to prove that $\VVV(T, (\bs{D}_{2}^{\infty})^{\ast})$ is a locally finite subvariety of $\VVV(\bs{RN})$, lacking the ES property. The fact that $\VVV(T, (\bs{D}_{2}^{\infty})^{\ast})$ is a subvariety of $\VVV(\bs{RN})$ follows from $T \subseteq \HHH\SSS(\bs{RN})$ and Lemma \ref{Lem:D2-in-RN}. 
Next we turn our attention to proving that
$\VVV(T, (\bs{D}_{2}^{\infty})^{\ast})$ is locally finite. 
Let $\D$ be the Heyting algebra depicted below:
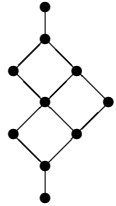
\begin{center}
\begin{picture}(44,80)(-4,-16)
\put(12,0){\line(0,-1){12}}
\put(12,-12){\circle*{4}}

\put(0,12){\line(1,-1){12}}
\put(0,12){\circle*{4}}
\put(0,12){\line(1,1){24}}
\put(12,0){\circle*{4}}
\put(12,0){\line(1,1){24}}
\put(24,12){\circle*{4}}
\put(24,12){\line(-1,1){24}}
\put(36,24){\circle*{4}}
\put(36,24){\line(-1,1){24}}
\put(12,24){\circle*{4}}
\put(24,36){\circle*{4}}
\put(0,36){\line(1,1){12}}
\put(0,36){\circle*{4}}
\put(12,48){\circle*{4}}

\put(12,48){\line(0,1){12}}
\put(12,60){\circle*{4}}
\end{picture}
\end{center}
Observe that the equation 
\[
 \bigvee_{i= 1}^{3} (x \to y_{i}) \lor (y_{i} \to x) \thickapprox 1
\]
holds in $\VVV(T, (\bs{D}_{2}^{\infty})^{\ast})$ but fails in $\D$. As a consequence, $\D \notin \VVV(T, (\bs{D}_{2}^{\infty})^{\ast})$. Since $\D$ has the form of one of the algebras in condition (ii) of the statement of Theorem \ref{Thm:locally-finite-KG}, we conclude that $\VVV(T, (\bs{D}_{2}^{\infty})^{\ast})$ is locally finite. Finally, we know that $T \cup \{ (\bs{D}_{2}^{\infty})^{\ast} \} \subseteq \class{W}_{2} \cap \class{ID}_{2}$. Together with Lemma \ref{Lem:excluding-D2}, this implies that $\VVV(T, (\bs{D}_{2}^{\infty})^{\ast})$ lacks the ES property.
\end{proof}

We conclude this section by showing that the ES property has an interesting consequence for subvarieties of $\class{KG}$:

\begin{Theorem}
\label{Thm:ES-implies-local-finiteness}
Let $\class{K}$ be a subvariety of $\class{KG}$. If $\class{K}$ has the ES property, then it is locally finite.
\end{Theorem}


\begin{proof}
We reason by contraposition. Suppose that a subvariety $\class{K}$ of $\class{KG}$ is not locally finite. From condition (iii) of Theorem \ref{Thm:locally-finite-KG} it follows that $\class{K}$ contains the sum $\bs{RN} + \bs{2}$. Now, it is easy to see that finitely generated subalgebras of $(\bs{D}_{2}^{\infty})^{\ast}$ belong to $\SSS\HHH(\bs{RN} + \bs{2})$. As a consequence, $(\bs{D}_{2}^{\infty})^{\ast} \in \III\SSS\PPU\SSS\HHH(\bs{RN} + \bs{2}) \subseteq \class{K}$. From Lemma \ref{Lem:excluding-D2}, it follows that $\class{K}$ lacks the ES property.
\end{proof}
It is show in \cite{BBdeJ08} 
that the cardinality of the interval $[\VVV(\bs{RN}),\class{KG}]$ is $2^{\aleph_0}$. No variety in this interval is locally finite, because they all contain $\bs{RN}$, which is an infinite finitely generated algebra. So, by Theorem~\ref{Thm:ES-implies-local-finiteness},
there is also a continuum of non-locally finite subvarieties of $\class{KG}$ without the ES property. 

\section{Epimorphism surjectivity in subvarieties of $\class{KG}$}

In obtaining a characterization of the subvarieties of $\class{KG}$ with the ES property, we shall rely on a series of technical observations. An element $a$ of a Heyting algebra $\A$ is said to be a \emph{node} if it is comparable with every element of $A$, see for instance \cite{Cit}.

\begin{Lemma}
\label{lem:KG}
Let $\A$ be a finite Heyting algebra and $0 < n \in \omega$.
\benroman
\item Suppose that $\A$ is one-generated. If $\vert A \vert \geq 6n + 1$, then
\[
\bs{2} + \underbrace{{\bs{D}_{2}}^{\ast} + \dots + {\bs{D}_{2}}^{\ast}}_{n\text{-times}} \in \SSS(\A).
\]
\item Suppose that $\A \in \class{KG}_{FSI}$, and that $c_{n+2} < c_{n+1} < \dots < c_{1}$ is a chain in $\A$, where $c_{n+2}, c_{n+1}, \dots, c_{1}$ are exactly the nodes of $\A$ in the interval $[c_{n+2}, c_{1}]$. If each interval $[c_{i+1}, c_{i}]$ has at least $3$ elements, then there are $\C_{1}, \dots, \C_{n} \in \{ {\bs{X}_{2}}^{\ast}, {\bs{D}_{2}}^{\ast} \}$%
\footnote{
Recall that the Esakia space $\bs{X}_2$ was defined before Remark~\ref{Rem:duals-of-Xn} on page~\pageref{Rem:duals-of-Xn}.
}
such that 
\[
\C_{1} + \dots + \C_{n} + \bs{2} \in \VVV(\A).
\]
\eroman
\end{Lemma}

\begin{proof}
(i): Recall that the class of one-generated Heyting algebras is $\HHH(\bs{RN})$. Hence $\A$ is a finite homomorphic image of $\bs{RN}$. Bearing this in mind, it is not hard to see that the poset reduct of $\A$ is isomorphic to a finite principal downset ${\downarrow} b$ of $\bs{RN}$. Accordingly, in what follows we shall identify the universes of $\A$ and ${\downarrow} b$, thus labelling the elements of $\A$ as in Figure~\ref{fig:RN} on page~\pageref{fig:RN}. A figure of ${\downarrow} b$ is reproduced below for convenience. Now, we define
\[
B \coloneqq \{ b, 0 \} \cup \bigcup_{n > k \in \omega} \{ w_{1 + 3k}, w_{2+3k}, a_{1+3k}\}.
\]
The elements of $B$ are encircled in the figure below.
Bearing in mind that $A = {\downarrow b}$ and $\vert A \vert \geq 6n +1$, it is easy to see that $B$ is the universe of a subalgebra of $\A$ isomorphic to $\bs{2} + \underbrace{{\bs{D}_{2}}^{\ast} + \dots + {\bs{D}_{2}}^{\ast}}_{n\text{-times}}$.

\begin{center}
\begin{picture}(62,133)
%
%
\put(38,3){\circle*{4}}
\put(38,3){\circle{8}}
\put(38,3){\line(1,1){12}}
\put(38,3){\line(-1,1){12}}
\put(50,15){\circle*{4}}
\put(50,15){\circle{8}}
\put(50,15){\line(-1,1){12}}
\put(26,15){\circle*{4}}
\put(26,15){\line(1,1){12}}
\put(38,27){\circle*{4}}
\put(12,5){$w_0$}
\put(50,6){$w_1$}
\put(42,25){$a_1$}

\put(26,15){\line(-1,1){12}}
\put(38,27){\line(-1,1){12}}
\put(14,27){\circle*{4}}
\put(14,27){\circle{8}}
\put(14,27){\line(1,1){12}}
\put(26,39){\circle*{4}}
\put(26,39){\circle{8}}
\put(0,17){$w_2$}
\put(30,37){$a_2$}

\put(38,27){\line(1,1){12}}
\put(50,39){\circle*{4}}
\put(50,39){\line(-1,1){12}}
\put(26,39){\line(1,1){12}}
\put(38,51){\circle*{4}}
\put(50,30){$w_3$}
\put(42,49){$a_3$}

\put(26,39){\line(-1,1){12}}
\put(38,51){\line(-1,1){12}}
\put(14,51){\circle*{4}}
\put(14,51){\circle{8}}
\put(14,51){\line(1,1){12}}
\put(26,63){\circle*{4}}
\put(0,41){$w_4$}
\put(30,61){$a_4$}

\put(38,51){\line(1,1){12}}
\put(50,63){\circle*{4}}
\put(50,63){\circle{8}}
\put(50,63){\line(-1,1){12}}
\put(26,63){\line(1,1){12}}
\put(38,75){\circle*{4}}
\put(38,75){\circle{8}}
\put(50,54){$w_5$}
\put(42,73){$a_5$}

\put(26,63){\line(-1,1){12}}
\put(38,75){\line(-1,1){12}}
\put(14,75){\circle*{4}}
\put(14,75){\line(1,1){12}}
\put(26,87){\circle*{4}}
\put(0,65){$w_6$}
\put(30,85){$a_6$}

\put(38,75){\line(1,1){12}}
\put(50,87){\circle*{4}}
\put(50,87){\circle{8}}
\put(50,87){\line(-1,1){10}}
\put(26,87){\line(1,1){10}}
\put(50,78){$w_7$}

\put(26,87){\line(-1,1){10}}

\put(32,105){\circle*{2}}
\put(32,111){\circle*{2}}
\put(32,117){\circle*{2}}

\put(32,129){\circle{8}}
\put(32,129){\circle*{4}}
\put(36,127){$b$}
\end{picture}

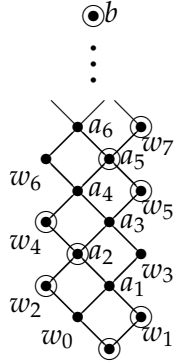
\captionof{figure}{$\bs{A}$ with the elements of $B$ encircled}
\label{fig:RN-subalgebra}
\end{center}

(ii): Since $\A$ is a finite FSI member of $\class{KG}$, we can apply Lemma \ref{Lem:finite-FSI-KG} obtaining that $\A = \B_{1} + \dots +  \B_{k}$ for some finite $\B_{1}, \dots, \B_{k} \in \HHH(\bs{RN})$. We can choose this decomposition of $\A$ into a sum in such a way that the unique nodes of $\A$ are 
\[
0^{\B_{k}} < 1^{\B_{k}} < \dots < 1^{\B_{1}}.
\]
Then $c_{n+2} < c_{n+1} < \dots < c_{1}$ is a segment of the above chain, so
\[
c_{1} = 1^{\B_{m}}, c_{2} = 1^{\B_{m+1}}, \dots, c_{n+1} = 1^{\B_{m +n}},
\]
for some $m \in \omega$ such that $m+n \leq k$, and $c_{n+2} = 0$ if $m+n=k$ and $c_{n+2} = 1^{\B_{m +n+1}}$ otherwise.

Then consider $m \leq i \leq m +n$. We have that $\B_{i}$ is a finite one-generated algebra. Moreover, the assumptions on the elements $c_{j}$, show that $\B_{i}$ has at least three elements and its unique nodes are the maximum and minimum elements. Now, the fact that $\B_{i}$ is one-generated implies that its underlying poset is isomorphic to a finite principal downset ${\downarrow} b$ of $\bs{RN}$, and in what follows we shall identify the universes of $\B_{i}$ and ${\downarrow} b$.

Since $\B_{i}$ has at least three elements and only two nodes, we obtain that, in the notation of Figure~\ref{fig:RN} on page~\pageref{fig:RN-subalgebra}, $b = a_{3p+q}$ for some $p \in \omega$ and $q \in \{ 0, 1, 2 \}$. Then for every $k \in \omega$ we define the following subsets of $RN$:
\begin{align*}
W \coloneqq & \{ w_{t} \colon t \in \omega \text{ is not divisible by }3 \}\\
C_{3k + 2} \coloneqq & \{ w_{t} \in W \colon t \leq 2 + 3k \} \cup \{ a_{3t + 2} \colon t \in \omega \text{ and } t \leq k \} \cup \{ 0, 1\}\\
C_{3k+1} \coloneqq & C_{3(k-1)+2} \cup \{ a_{1 + 3k}, a_{3k}, w_{1 + 3k} \}\\
C_{3k} \coloneqq & C_{3(k-1)+2} \cup \{ a_{3k}, w_{3k}, a_{3k-2}, a_{3(k-1)} \}. 
\end{align*}
We have that $C_{3p+q}$ is the universe of a subalgebra $\C_{i}$ of $\B_{i}$, which is isomorphic to a sum of the form $\C^{i}_{1} + \dots + \C^{i}_{k_{i}}$ for some $\C^{i}_{1}, \dots, \C^{i}_{k_{i}} \in \{ {\bs{X}_{2}}^{\ast}, {\bs{D}_{2}}^{\ast} \}$. See below for figures illustrating $\B_i$ for the three possibilities of $q$, where the elements of $C_{3p+q}$ are encircled.

This easily implies that
\[
\C_{m} + \dots + \C_{m + (n-1)} + \bs{2} \in \SSS(\C_{m} + \dots + \C_{m+n}) \subseteq \SSS(\B_{m} + \dots + \B_{m+n}).
\]
Bearing in mind that $\A = \B_{1} + \dots +  \B_{k}$, this yields that
\[
\C_{m} + \dots + \C_{m + (n-1)} + \bs{2} \in \SSS\HHH(\A) \subseteq \VVV(\A).
\]
Putting all this together we get that
\[
\C^{m}_{1}+ \dots + \C^{m}_{k_{m}}+\C^{m+1}_1 + \dots + \C^{m+(n-1)}_{k_{m+(n-1)}}+\bs{2} \in \VVV(\A).
\]
Finally, we let $\D_{1}, \dots, \D_{n} \in \{ {\bs{X}_{2}}^{\ast}, {\bs{D}_{2}}^{\ast} \}$ be the first $n$ components of the sum above, so that
\[
\D_{1} + \dots + \D_{n} + \bs{2} \in \SSS(\C_{m} + \dots + \C_{m + (n-1)} + \bs{2}) \subseteq \VVV(\A),
\]
concluding the proof.
\end{proof}

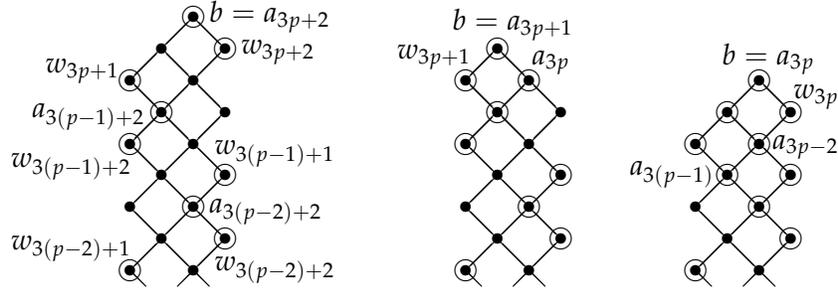
\begin{center}
\begin{tabular}{ccccc}

\begin{picture}(126,113)(0,20)
%
%


\qbezier(69,29)(69,29)(75,23)
\qbezier(69,29)(69,29)(63,23)
\put(69,29){\circle*{4}}

\put(0,37){$w_{3(p-2)+1}$}
\put(45,29){\circle*{4}}
\put(45,29){\circle{8}}
\put(45,29){\line(1,1){12}}
\qbezier(45,29)(45,29)(51,23)
\put(57,41){\circle*{4}}
\put(57,41){\line(1,-1){12}}

\put(57,41){\line(1,1){12}}
\put(77,29){$w_{3(p-2)+2}$}
\put(81,41){\circle*{4}}
\put(81,41){\circle{8}}
\put(81,41){\line(-1,1){12}}
\put(81,41){\line(-1,-1){12}}
\put(75,51){$a_{3(p-2)+2}$}
\put(69,53){\circle*{4}}
\put(69,53){\circle{8}}

\put(45,53){\circle*{4}}
\put(45,53){\line(1,1){12}}
\put(45,53){\line(1,-1){12}}
\put(57,65){\circle*{4}}
\put(57,65){\line(1,-1){12}}

\put(57,65){\line(1,1){12}}
\put(77,73){$w_{3(p-1)+1}$}
\put(81,65){\circle*{4}}
\put(81,65){\circle{8}}
\put(81,65){\line(-1,1){12}}
\put(81,65){\line(-1,-1){12}}
\put(69,77){\circle*{4}}

\put(0,68){$w_{3(p-1)+2}$}
\put(45,77){\circle*{4}}
\put(45,77){\circle{8}}
\put(45,77){\line(1,1){12}}
\put(45,77){\line(1,-1){12}}
\put(8,87){$a_{3(p-1)+2}$}
\put(57,89){\circle*{4}}
\put(57,89){\circle{8}}
\put(57,89){\line(1,-1){12}}

\put(57,89){\line(1,1){12}}
\put(81,89){\circle*{4}}
\put(81,89){\line(-1,1){12}}
\put(81,89){\line(-1,-1){12}}
\put(69,101){\circle*{4}}

\put(13,104){$w_{3p+1}$}
\put(45,101){\circle*{4}}
\put(45,101){\circle{8}}
\put(45,101){\line(1,1){12}}
\put(45,101){\line(1,-1){12}}
\put(57,113){\circle*{4}}
\put(57,113){\line(1,-1){12}}

\put(57,113){\line(1,1){12}}
\put(87,111){$w_{3p+2}$}
\put(81,113){\circle*{4}}
\put(81,113){\circle{8}}
\put(81,113){\line(-1,1){12}}
\put(81,113){\line(-1,-1){12}}
\put(75,123){$b=a_{3p+2}$}
\put(69,125){\circle*{4}}
\put(69,125){\circle{8}}

%
%
%
%

\end{picture}

&&

\begin{picture}(67,110)(19,20)
%
%


\qbezier(69,29)(69,29)(75,23)
\qbezier(69,29)(69,29)(63,23)
\put(69,29){\circle*{4}}

\put(45,29){\circle*{4}}
\put(45,29){\circle{8}}
\put(45,29){\line(1,1){12}}
\qbezier(45,29)(45,29)(51,23)
\put(57,41){\circle*{4}}
\put(57,41){\line(1,-1){12}}

\put(57,41){\line(1,1){12}}
\put(81,41){\circle*{4}}
\put(81,41){\circle{8}}
\put(81,41){\line(-1,1){12}}
\put(81,41){\line(-1,-1){12}}
\put(69,53){\circle*{4}}
\put(69,53){\circle{8}}

\put(45,53){\circle*{4}}
\put(45,53){\line(1,1){12}}
\put(45,53){\line(1,-1){12}}
\put(57,65){\circle*{4}}
\put(57,65){\line(1,-1){12}}

\put(57,65){\line(1,1){12}}
\put(81,65){\circle*{4}}
\put(81,65){\circle{8}}
\put(81,65){\line(-1,1){12}}
\put(81,65){\line(-1,-1){12}}
\put(69,77){\circle*{4}}

\put(45,77){\circle*{4}}
\put(45,77){\circle{8}}
\put(45,77){\line(1,1){12}}
\put(45,77){\line(1,-1){12}}
\put(57,89){\circle*{4}}
\put(57,89){\circle{8}}
\put(57,89){\line(1,-1){12}}

\put(57,89){\line(1,1){12}}
\put(81,89){\circle*{4}}
\put(81,89){\line(-1,1){12}}
\put(81,89){\line(-1,-1){12}}
\put(70,107){$a_{3p}$}
\put(69,101){\circle*{4}}
\put(69,101){\circle{8}}

\put(19,109){$w_{3p+1}$}
\put(45,101){\circle*{4}}
\put(45,101){\circle{8}}
\put(45,101){\line(1,1){12}}
\put(45,101){\line(1,-1){12}}
\put(40,120){$b=a_{3p+1}$}
\put(57,113){\circle*{4}}
\put(57,113){\circle{8}}
\put(57,113){\line(1,-1){12}}

\end{picture}

&&

\begin{picture}(81,98)(19,20)
%
%


\qbezier(69,29)(69,29)(75,23)
\qbezier(69,29)(69,29)(63,23)
\put(69,29){\circle*{4}}

\put(45,29){\circle*{4}}
\put(45,29){\circle{8}}
\put(45,29){\line(1,1){12}}
\qbezier(45,29)(45,29)(51,23)
\put(57,41){\circle*{4}}
\put(57,41){\line(1,-1){12}}

\put(57,41){\line(1,1){12}}
\put(81,41){\circle*{4}}
\put(81,41){\circle{8}}
\put(81,41){\line(-1,1){12}}
\put(81,41){\line(-1,-1){12}}
\put(69,53){\circle*{4}}
\put(69,53){\circle{8}}

\put(45,53){\circle*{4}}
\put(45,53){\line(1,1){12}}
\put(45,53){\line(1,-1){12}}
\put(20,65){$a_{3(p-1)}$}
\put(57,65){\circle*{4}}
\put(57,65){\circle{8}}
\put(57,65){\line(1,-1){12}}

\put(57,65){\line(1,1){12}}
\put(81,65){\circle*{4}}
\put(81,65){\circle{8}}
\put(81,65){\line(-1,1){12}}
\put(81,65){\line(-1,-1){12}}
\put(74,75){$a_{3p-2}$}
\put(69,77){\circle*{4}}
\put(69,77){\circle{8}}

\put(45,77){\circle*{4}}
\put(45,77){\circle{8}}
\put(45,77){\line(1,1){12}}
\put(45,77){\line(1,-1){12}}
\put(57,89){\circle*{4}}
\put(57,89){\circle{8}}
\put(57,89){\line(1,-1){12}}

\put(57,89){\line(1,1){12}}
\put(82,94){$w_{3p}$}
\put(81,89){\circle*{4}}
\put(81,89){\circle{8}}
\put(81,89){\line(-1,1){12}}
\put(81,89){\line(-1,-1){12}}
\put(55,107){$b = a_{3p}$}
\put(69,101){\circle*{4}}
\put(69,101){\circle{8}}

\end{picture}

\end{tabular}
\end{center}

\begin{Lemma}\label{Lem:(iii)-->(i)}.
Let $\mathsf{K}$ be a subvariety of $\mathsf{KG}$. Suppose there exists a positive integer $n$ such that $\mathsf{K}$ excludes all sums of the form $\A_1 + \dots + \A_{n} + \bs{2}$, where $\A_i \in \{{\D_2}^{\ast}, {\bs{X}_2}^{\ast}\}$ for every $i \leq n$. Then $\mathsf{K}$ has the ES property.
\end{Lemma}

\begin{proof}

We claim that $\class{K}$ is locally finite. Suppose the contrary. By Theorem \ref{Thm:locally-finite-KG}, we obtain that $\bs{RN} + \bs{2} \in \class{K}$. Now, it is easy to see that
\[
\bs{2}+  \underbrace{{\bs{D}_{2}}^{\ast} + \dots + {\bs{D}_{2}}^{\ast}}_{n\text{-times}} + \bs{2} \in \SSS(\bs{RN} + \bs{2})
\]
and, therefore, 
\[
\underbrace{{\bs{D}_{2}}^{\ast} + \dots + {\bs{D}_{2}}^{\ast}}_{n\text{-times}} + \bs{2} \in \HHH(\bs{2}+  \underbrace{{\bs{D}_{2}}^{\ast} + \dots + {\bs{D}_{2}}^{\ast}}_{n\text{-times}} + \bs{2}) \subseteq \SSS\HHH(\bs{RN} + \bs{2}) \subseteq \class{K}.
\]
But this contradicts the assumption, thus establishing the claim.

Now, suppose with a view to obtaining a contradiction 
that $\mathsf{K}$ does not have the ES property. Then by Theorem \ref{Thm:Campercholi} there exists an FSI algebra $\B \in \mathsf{K}$ with a proper $\mathsf{K}$-epic subalgebra $\A$. Clearly, $\B$ embeds into an ultraproduct $\C = \prod_{i \in I} \C_i/\mathcal{U}$ of its finitely generated subalgebras. Henceforth, we identify $\B$ with its image under this embedding (so that both $\A$ and $\B$ are subalgebras of $\C$). Observe that, since $\class{K}$ is locally finite, we know that each $\C_i$ is finite. Moreover, by Lemma \ref{Lem:correspondences}(\ref{Lem:correspondences:FSI}), each $\C_i$ is FSI.  By Lemma \ref{Lem:finite-FSI-KG} each $\C_i$ is a finite sum of finite algebras in $\HHH(\bs{RN})$.

\textbf{Claim~(a).} If $c_1 < \dots < c_{n+2}$ is a sequence of nodes in $\C$ such that $c_1, \dots, c_{n+2}$ are exactly the nodes of $\C$ in the interval $[c_{1}, c_{n+2}]$, then there is at least one $j \in \{ 1, \dots, n+1\}$ such that the interval $[c_j, c_{j+1}]$ is a two-element set.

\emph{Proof of Claim~(a).} As the statement of Claim~(a) can be formulated as a first-order sentence in the language of Heyting algebras, it follows by \L o\'s' Theorem \cite[Sec.\ V, Thm.\ 2.9]{BS81} that Claim~(a) will hold in $\C$ if it holds in $\C_i$ for every $i \in I$. Suppose on the contrary that there is an $i \in I$ and a sequence of nodes $a_1 < \dots < a_{n+2}$ of $\C_{i}$ as above such that each interval $[a_{j}, a_{j+1}]$ has at least three elements. Together with Lemma \ref{lem:KG}(ii), this implies that there are $\bs{B}_{1}, \dots, \bs{B}_{n} \in \{ {\bs{X}_{2}}^{\ast}, {\bs{D}_{2}}^{\ast} \}$ such that
\[
\bs{B}_{1} + \dots + \bs{B}_{n} + \bs{2} \in \VVV(\C_{i}) \subseteq \class{K}.
\]
But this contradicts the assumptions and, therefore, establishes the claim.

\textbf{Claim (b).} For every pair of nodes $c < d$ in $\C$ with no other node in between, we have that $\vert [c,d] \vert \leq 6n+6$.

\emph{Proof of Claim (b).} As with Claim~(a), it suffices to show that each $\C_i$ satisfies this property. To this end, consider nodes $a < b$ of $\C_{i}$ with no other node in between, and suppose, with a view to obtaining a contradiction, 
that $\vert [a, b] \vert \geq 6 + 7$. By Lemma \ref{Lem:finite-FSI-KG}, we have that $\C_{i} = \B_{1} +  \dots + \B_{m}$ for some $\B_{1}, \dots, \B_{m} \in \HHH(\bs{RN})$. Since the unique nodes in the interval $[a, b]$ are $a$ and $b$, there exists $j \leq m$ such that $[a, b] \subseteq \B_{j}$. Hence $\B_{j}$ is a one-generated Heyting algebra such that $\vert B \vert \geq 6(n+1) + 1$. By Lemma \ref{lem:KG}(i) we obtain that
\[
\bs{2} + \underbrace{{\bs{D}_{2}}^{\ast} + \dots + {\bs{D}_{2}}^{\ast}}_{n+1\text{-times}} \in \SSS(\B_{j}).
\]
This easily implies that
\[
\bs{2} + \underbrace{{\bs{D}_{2}}^{\ast} + \dots + {\bs{D}_{2}}^{\ast}}_{n\text{-times}} + \bs{2} \in \SSS(\C_{i}).
\]
As a consequence, we obtain
\[
\underbrace{{\bs{D}_{2}}^{\ast} + \dots + {\bs{D}_{2}}^{\ast}}_{n\text{-times}} + \bs{2}  \in \HHH(\bs{2} + \underbrace{{\bs{D}_{2}}^{\ast} + \dots + {\bs{D}_{2}}^{\ast}}_{n\text{-times}} + \bs{2} ) \subseteq \HHH\SSS(\C_{i}) \subseteq \class{K}.
\]
But this contradicts the assumptions and, therefore, establishes the claim.

\textbf{Claim (c).} For every element $c$ of $\C$ there exist a largest node of $\C$ below $c$, and a smallest node of $\C$ above $c$.

\emph{Proof of Claim (c).} As with Claim~(a), it suffices to show that each $\C_i$ satisfies this property. But this is an immediate consequence of the fact that $\C_{i}$ is finite, and we are done.

Now, we will extend $\A$ to a proper subalgebra $\D \leq \B$ such that the inclusion map $\D \to \B$ is almost onto. This will contradict the fact that $\class{K}$ has the weak ES property (see Theorem \ref{Thm:Kreisel}). In order to construct the extension $\D$ of $\A$, we reason as follows. Since $\A$ is a proper subalgebra of $\B$, we can choose an element $b \in B \smallsetminus A$. By Claim~(c), there exist a largest node $c$ of $\C$ such that $c \leq b$ and a smallest node $d$ of $\C$ such that $b \leq d$.

If ${\uparrow} d$ is finite, we set $\A' \coloneqq \A$. Now assume that ${\uparrow}d$ is infinite. We show, by supposing the contrary, 
that there exists a sequence of nodes $d = b_1 < b_2 < \dots < b_{n+2}$, where $b_1, \dots, b_{n+2}$ are exactly the nodes of $\A$ in the interval $[b_1,b_{n+2}]$. Then there are only finitely many nodes above $b$. By Claim~(c) every element of ${\uparrow} b$ belongs to an interval between nodes, and there are only finitely many such intervals. But by Claim~(b) this implies that ${\uparrow}b$ is finite, which contradicts the assumption, thus establishing the existence of the sequence. Then, by Claim~(a), there exists $j \leq n$ such that $[b_j,b_{j+1}]$ is a two-element set. Let $A' = A \cup {\uparrow} b_{j+1}$. Since $b_{j+1}$ is a node, it is not hard to see that $A'$ is the universe of a subalgebra $\A'$ of $\B$. Moreover, by Claim~(b), the interval $[b,b_{j+1}]$ is finite.

Now, if ${\downarrow}c$ is finite we let $\D \coloneqq \A'$. Then suppose that ${\downarrow}c$ is infinite. Using an argument similar to the one above, we can construct a chain of nodes $c = c_1 > c_2 > \dots > c_{n+2}$, where $c_1, \dots, c_{n+2}$ are exactly the nodes of $\A'$ in the interval $[c_{n+2},c_1]$. By Claim~(a), there exists $j \leq n+1$ such that $[c_{j+1},c_j]$ is a two-element set. Again the set $D = A' \cup {\downarrow}c_{j+1}$ is the universe of a subalgebra $\D$ of $\B$. Moreover, $[c_{j+1},c]$ is finite by Claim~(b).

Observe that the subalgebra $\D \leq \B$ extends $\A$ and is proper, since $b \notin D$. Moreover, $\D$ is a $\class{K}$-epic subalgebra of $\B$, since $\D$ extends $\A$. Bearing in mind that the interval $[c, d]$ is finite by Claim~(b), we obtain that $B \smallsetminus D$ is finite. Hence the inclusion map $\D \to \B$ is an almost onto non-surjective $\class{K}$-epimorphism. This implies that $\class{K}$ lacks the weak ES property, contradicting Theorem \ref{Thm:Kreisel}.
\end{proof}

\begin{Lemma}\label{Lem:(i)-->(ii)}
Let $\{ \bs{Z}_{n} \colon n \in \omega \}$ be a family of Esakia spaces such that $\bs{Z}_{n} \in \{ \bs{X}_{2}, \bs{D}_{2} \}$ for every $n \in \omega$, and let $\class{K}$ be a subvariety of $\class{ID}_{2} \cap \class{W}_{2}$. If $\sum \bs{Z}_{n} \in \class{K}_{\ast}$, then $\class{K}$ lacks the ES property.
\end{Lemma}

\begin{proof}
Suppose with a view to obtaining a contradiction %
that there is a subvariety $\class{K}$ of $\class{ID}_{2} \cap \class{W}_{2}$ with the ES property and such that $\sum \bs{Z}_{n} \in \class{K}_{\ast}$ for some family $\{ \bs{Z}_{n} \colon n \in \omega \}$ of Esakia spaces as in the statement.

We claim that for every $m \in \omega$, there are $t, k \geq m$ such that $\bs{Z}_{t} = \bs{D}_{2}$ and $\bs{Z}_{k} = \bs{X}_{2}$. Suppose not. Then there exists $m \in \omega$ such that either $\bs{Z}_{k} = \bs{D}_{2}$ for every $k \geq m$ or $\bs{Z}_{k} = \bs{X}_{2}$ for every $k \geq m$. This implies that either $\bs{D}_{2}^{\infty}$ or $\bs{X}_{2}^{\infty}$ is an E-subspace of $\sum \bs{Z}_{n}$ and, therefore, that either $\bs{D}_{2}^{\infty} \in \class{K}_{\ast}$ or $\bs{X}_{2}^{\infty} \in \class{K}_{\ast}$. By Theorem \ref{Thm:Es-bounded-width} and Lemma \ref{Lem:excluding-D2} we conclude that $\class{K}$ lacks the ES property, thus contradicting the assumptions. This establishes the claim.

We know that every component $\bs{Z}_{m}$ of $\sum \bs{Z}_{n}$ is a copy of either $\bs{D}_{2}$ or $\bs{X}_{2}$. Accordingly, we denote the elements of $\bs{Z}_{m}$ as follows:
\begin{center}
\begin{tabular}{ccc}
\begin{picture}(50,45)(-5,-20)
\put(16,8){\line(1,1){12}}
\put(16,20){\line(1,-1){12}}
\put(28,8){\line(0,1){12}}
\put(16,8){\circle*{4}}
\put(28,8){\circle*{4}}
\put(16,20){\circle*{4}}
\put(28,20){\circle*{4}}
\put(0,4){$a_m$}
\put(32,4){$b_m$}
\put(0,18){$c_m$}
\put(32,18){$d_m$}

\put(0,-2){\rotatebox{270}{\scalebox{1}[4]{$\}$}}}
\put(-4,-18){if $\bs{Y}_m = \bs{X}_2$}
\end{picture}
&\hspace{\textheight}&
\begin{picture}(50,45)(-5,-20)
\put(16,8){\circle*{4}}
\put(28,8){\circle*{4}}
\put(0,4){$a_m$}
\put(32,4){$b_m$}

\put(0,-2){\rotatebox{270}{\scalebox{1}[4]{$\}$}}}
\put(-4,-18){if $\bs{Y}_m = \bs{D}_2$}
\end{picture}
\end{tabular}
\end{center}
Then consider the relation $R$ on $\sum Z_{n}$ defined as follows: for every $x, y \in \sum Z_{n}$,
\begin{align*}
\langle x, y \rangle \in R \Longleftrightarrow & \text{ either }x = y\\
& \text{ or }(x, y \in \bs{Z}_{m} = \bs{X}_{2} \text{ and }\{ x, y \} =\{ b_{m}, d_{m}\}, \text{ for some $m \in \omega$})\\
& \text{ or }(x \in \bs{Z}_{m} = \bs{X}_{2}, y \in \bs{Z}_{m+1} = \bs{D}_{2} \text{ and }\{ x, y \} =\{ c_{m}, a_{m+1}\},\\
&\ \ \ \ \ \text{ for some $m \in \omega$})\\
& \text{ or }(x \in \bs{Z}_{m} = \bs{D}_{2}, y \in \bs{Z}_{m+1} = \bs{X}_{2} \text{ and }\{ x, y \} =\{ b_{m}, a_{m+1}\},\\
&\ \ \ \ \ \text{ for some $m \in \omega$})\\
& \text{ or }(x \in \bs{Z}_{m} = \bs{D}_{2}, y \in \bs{Z}_{m+1} = \bs{D}_{2} \text{ and }\{ x, y \} =\{ b_{m}, a_{m+1}\},\\
&\ \ \ \ \ \text{ for some $m \in \omega$})\\
& \text{ or }(x \in \bs{Z}_{m} = \bs{X}_{2}, y \in \bs{Z}_{m+1} = \bs{X}_{2} \text{ and }\{ x, y \} =\{ c_{m}, a_{m+1}\},\\
&\ \ \ \ \ \text{ for some $m \in \omega$}).
\end{align*}
An argument, similar to the one detailed for Lemma \ref{Lem:correct-partition}, shows that $R$ is a correct partition on $\sum \bs{Z}_{n}$.

Since $\class{K}$ has the ES property, we can apply Lemma \ref{Lem:epic-subalgebras-dual}, obtaining that there exist $\boldsymbol{Y} \in \class{K}_{\ast}$ and a pair of different Esakia morphisms $f, g \colon \boldsymbol{Y} \to \sum \bs{Z}_{n}$ such that  $\langle f(y), g(y) \rangle \in R$ for every $y \in Y$. Since $f \ne g$, there exists $\bot \in Y$ such that $f(\bot) \ne g(\bot)$. We can assume without loss of generality that $f(\bot) < g(\bot)$. Moreover, as in the proof of Theorem \ref{Thm:Es-bounded-width}, we can assume that $Y = {\uparrow} \bot$. Observe that $\sum \bs{Z}_{n}$ has width at most $2$, as does $\bs{Y}$, since $\class{K} \subseteq \class{W}_2$. Moreover, $\bs{Y}$ and $\sum \bs{Z}_{n}$, and the Esakia morphism $f \colon \boldsymbol{Y} \to \sum \bs{Z}_{n}$ satisfy the assumptions of Lemma \ref{Lem:trick-width}. Therefore, $\boldsymbol{Y}$ has a subposet $\langle Z; \leq^{\boldsymbol{Y}}\rangle$ such that the restriction
\[
f \colon \langle Z; \leq^{\boldsymbol{Y}}\rangle \to \langle {\uparrow} f(\bot)^{\top}; \leq^{\sum \bs{Z}_{n}}\rangle
\]
is a poset isomorphism. For the sake of simplicity, we denote the elements of $Z$ exactly as their alter egos in ${\uparrow} f(\bot)^{\top}$. It is not hard to see that
\[
f(z) = z \text{ and }g(z) = \max(z/R)
\]
for every $z \in Z$ (see the proof of Theorem \ref{Thm:Es-bounded-width}, if necessary).

From the claim it follows that there exists $m \in \omega$ such that $Z_{m} \cup Z_{m+1} \subseteq Z$, $\bs{Z}_{m} = \bs{D}_{2}$ and $\bs{Z}_{m+1} = \bs{X}_{2}$. The following picture represents the relevant part of $\sum \bs{Z}_{n}$ equipped with the correct partition $R$:
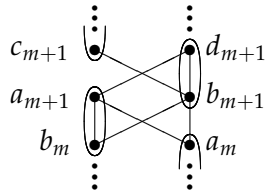
\begin{center}
\begin{picture}(100,85)(-5,-20)
\put(26,4){\line(2,1){36}}
\put(26,4){\circle*{4}}
\put(26,4){\line(0,1){18}}
\put(62,4){\circle*{4}}
\put(62,4){\line(0,1){18}}
\put(62,4){\line(-2,1){36}}

\put(26,22){\line(2,1){36}}
\put(26,22){\circle*{4}}
\put(62,22){\circle*{4}}
\put(62,22){\line(0,1){18}}
\put(62,22){\line(-2,1){36}}

\put(26,40){\circle*{4}}
\put(62,40){\circle*{4}}

\newsavebox{\EqRb}
\savebox{\EqRb}(8,26){
\qbezier(4,0)(0,0)(0,13)
\qbezier(0,13)(0,26)(4,26)
\qbezier(4,0)(8,0)(8,13)
\qbezier(8,13)(8,26)(4,26)%
}
\put(22,0){\usebox{\EqRb}}
\put(58,18){\usebox{\EqRb}}

\qbezier(26,36)(22,36)(22,49)
\qbezier(26,36)(30,36)(30,49)

\qbezier(58,-5)(58,8)(62,8)
\qbezier(66,-5)(66,8)(62,8)%

\put(26,48){\circle*{2}}
\put(26,52){\circle*{2}}
\put(26,56){\circle*{2}}

\put(26,-4){\circle*{2}}
\put(26,-8){\circle*{2}}
\put(26,-12){\circle*{2}}

\put(62,48){\circle*{2}}
\put(62,52){\circle*{2}}
\put(62,56){\circle*{2}}

\put(62,-4){\circle*{2}}
\put(62,-8){\circle*{2}}
\put(62,-12){\circle*{2}}

\put(5,2){$b_{m}$}
\put(-6,20){$a_{m+1}$}
\put(-6,38){$c_{m+1}$}
\put(68,2){$a_{m}$}
\put(68,20){$b_{m+1}$}
\put(68,38){$d_{m+1}$}
\end{picture}
\end{center}
Observe that $g(a_m) = a_{m} \leq^{\sum \bs{Z}_{n}}c_{m+1}$. Since $g$ is an Esakia morphism, there exists $y \in Y$ such that $a_{m} \leq^{\bs{Y}} y$ and $g(y) = c_{m+1}$. Moreover, $g(a_{m+1}) = a_{m+1}$ and $g(b_{m+1}) = d_{m+1}$. Together with the fact that $g(y) = c_{m+1}$, this implies that $g(y)$ is incomparable with $g(a_{m+1})$ and $g(b_{m+1})$ in $\sum \bs{Z}_{n}$. Since $g$ is order-preserving, we conclude that $y$ is incomparable with $a_{m+1}$ and $b_{m+1}$ in $\bs{Y}$. As $\bs{Y}$ has width $\leq 2$, we obtain that $a_{m+1}$ and $b_{m+1}$ are comparable is $\bs{Y}$, a contradiction.
\end{proof}

We are finally ready to present a characterization of the subvarieties of $\class{KG}$ with the ES property:

\begin{Theorem}\label{Thm:ES-in-KG}
Let  $\class{K} \subseteq \class{KG}$ be a variety. The following conditions are equivalent:
\benroman
\item $\class{K}$ has the ES property.
\item $\class{K}$ excludes all sums $\sum \A_{n}$ of families $\{ \A_{n} \colon n \in \omega \}$ of Heyting algebras such that $\A_n \in \{{\D_2}^{\ast}, {\bs{X}_2}^{\ast}\}$ for every $n \in \omega$.
\item There is a positive integer $n$ such that $\mathsf{K}$ excludes all sums of the form $\A_1 + \dots + \A_{n} + \bs{2}$, where $\A_i \in \{{\D_2}^{\ast}, {\bs{X}_2}^{\ast}\}$ for every $i \leq n$.
\eroman
\end{Theorem}

\begin{proof}
Lemmas \ref{Lem:(i)-->(ii)} and \ref{Lem:(iii)-->(i)} yield directions (i)$\Rightarrow$(ii) and (iii)$\Rightarrow$(i).

(ii)$\Rightarrow$(iii): Consider the expansion of the language of Heyting algebras with fresh constants $\{ a_{n} \colon n \in \omega \}$. Let $\class{node}(x)$ be the first-order formula whose meaning is ``$x$ is a node''. Similarly, let $\class{interval}(x, y)$ be the first-order formula meaning ``the interval $[x, y]$ is order-isomorphic to the poset reduct of either ${\D_{2}}^{\ast}$ or ${\bs{X}_{2}}^{\ast}$''. Consider the following set of sentences of the expanded language:
\[
\Phi \coloneqq \{ \class{node}(a_{n}) \colon n \in \omega \} \cup \{ \class{interval}(a_{n+1}, a_{n}) \colon n \in \omega \} \cup \{ a_{0} \thickapprox 1 \}.
\]
Moreover, let $\Delta$ the set of axioms of $\class{KG}$ (regarded as first-order sentences).

We shall reason by contraposition. Suppose that for every $n \in \omega$ there are $\A_1, \dots, \A_{n} \in \{{\D_2}^{\ast}, {\bs{X}_2}^{\ast}\}$ such that $\A_1 + \dots + \A_{n} + \bs{2} \in \class{K}$. This easily implies that every finite subset of $\Phi \cup \Delta$ has a model. By the Compactness Theorem, it follows that $\Phi \cup \Delta$ has a model $\A$. Clearly, the Heyting algebra reduct $\A^{-}$ of $\A$ belongs to $\class{K}$. Moreover, $\A^{-}$ has a subalgebra isomorphic to a sum $\sum \A_{n}$, where $\{ \A_{n} \colon n \in \omega \}$ is a family of Heyting algebras such that $\A_n \in \{{\D_2}^{\ast}, {\bs{X}_2}^{\ast}\}$ for every $n \in \omega$. 
\end{proof}

As a consequence we obtain an alternative proof of the following result from \cite[Cor.\ 5.7]{BMR17}:

\begin{Corollary}
Every variety of G\"odel algebras has the ES property.
\end{Corollary}

\begin{proof}
Let $\class{K}$ be a variety of G\"odel algebras. Observe that $\class{K}$ excludes every sum of the form $\A + \bs{B} + \bs{2}$ with $\A, \B \in \{ {\bs{X}_{2}}^{\ast}, {\bs{D}_{2}}^{\ast} \}$, since these sums are not contained in $\class{W}_1$. As $\class{K} \subseteq \class{KG}$, we can apply Theorem \ref{Thm:ES-in-KG}, obtaining that $\class{K}$ has the ES property. 
\end{proof}

In summary, we have shown, in Theorem~\ref{Thm:continuum}, that there is a continuum of locally finite subvarieties of $\VVV(\bs{RN}) \subseteq \class{KG}$ without the ES property. In Theorem~\ref{Thm:ES-in-KG}, we classified the subvarieties of $\class{KG}$ that have the ES property. In particular, they are all locally finite. Furthermore, notice that if a subvariety $\class{K}$ of $\class{KG}$ has the ES property, then so do all subvarieties of $\class{K}$.

\section{Brouwerian algebras}

Subreducts of Heyting algebras in the signature $\langle \land, \lor, \to, 1 \rangle$ are called \emph{Brouwerian algebras}. It is well-known that the varieties of Brouwerian algebras algebraize the positive logics, i.e., the axiomatic extensions of the $\langle \land, \lor, \to, 1 \rangle$-fragment of intuitionistic logic. As a consequence, a positive logic $\vdash$ has the infinite (deductive) Beth (definability) property if and only if the variety of Brouwerian algebras associated with $\vdash$ has the ES property.

This prompts the question of whether it is possible to adapt the results on epimorphisms obtained so far to varieties of Brouwerian algebras. In the majority of cases this can be done by changing naturally the notion of sums of algebras and the duality to the case of Brouwerian algebras 
(for the latter, see for instance \cite{BMR17}).

One exception is the proof of Theorem \ref{Thm:continuum}, which does not immediately survive the deletion of $0$ from the type. However, it remains true that there is a continuum of locally finite varieties of Brouwerian algebras lacking the ES property. The proof of this fact is an obvious adaptation of that of Theorem \ref{Thm:continuum}, in which the algebra $\B_{n}$ must be replaced (in the notation of the proof) by $\bs{2} + \A + \C_{1} + \dots +  \C_{n-2} + \A$.

\

\paragraph{\bfseries Acknowledgements.}
We thank Guram Bezhanishvili and James G.\ Raftery for helpful conversations on the topic. We are also very grateful to the anonymous referees, whose useful remarks helped to improve the presentation of the paper. The first author was supported by project CZ$.02.2.69/0.0/0.0/17\_050/0008361$, OPVVV M\v{S}MT, MSCA-IF Lidsk\'{e} zdroje v teoretick\'{e} informatice and by the Beatriz Galindo grant BEAGAL$18$/$00040$ of the Spanish Ministry of Science, Innovation and Universities.
The second author was supported by the DST-NRF Centre of Excellence in Mathematical and Statistical Sciences (CoE-MaSS), South Africa.
Opinions expressed and conclusions arrived at are those of the authors and are not necessarily to be attributed to the CoE-MaSS.

\bibliographystyle{plain}

\end{document}